\documentclass[oneside,10pt,a4paper,reqno]{amsart}
\usepackage{cd_comm}

  \let\origmaketitle\maketitle
\def\maketitle{
  \begingroup
  \def\uppercasenonmath##1{} 
  \let\MakeUppercase\relax 
  \origmaketitle
  \endgroup
}

\begin{document}

	\title{On Classical Determinate Truth}
\author{Luca Castaldo and Carlo Nicolai}
		\date{}
	
	\maketitle

 \begin{abstract}
The paper proposes and studies new classical, type-free theories of truth and determinateness with unprecedented features. The theories are fully compositional, strongly classical (namely, their internal and external logics are both classical), and feature a \emph{defined} determinateness predicate satisfying desirable and widely agreed principles. The theories capture a conception of truth and determinateness according to which the generalizing power associated with the classicality and full compositionality of truth is combined with the identification of a natural class of sentences -- the determinate ones -- for which clear-cut semantic rules are available. Our theories can also be seen as the \emph{classical closures} of Kripke-Feferman truth: their $\omega$-models, which we precisely pinned down, result from including in the extension of the truth predicate the sentences that are satisfied by a Kripkean closed-off fixed point model.
The theories compare to recent theories proposed by Fujimoto and Halbach, featuring a primitive determinateness predicate. In the paper we show that our theories entail all principles of Fujimoto and Halbach's theories, and are proof-theoretically equivalent to Fujimoto and Halbach's $\cdplus$. {We also show establish some negative results on Fujimoto and Halbach's theories: such results show that, unlike what happens in our theories, the primitive determinateness predicate prevents one from establishing clear and unrestricted semantic rules for the language with type-free truth. }
\end{abstract}

\section{A Conception of Truth (and Determinateness)}

\noindent In this work we offer a cluster of formal theories with unprecedented features and prove several results about them. The theories are accompanied by a conception of truth (and determinateness), which we outline in this introductory section.\footnote{A full exposition and philosophical defence of the conception is carried out in a companion paper.}

\subsection{Desiderata for truth}

Following a long-standing tradition in formal theories of truth, starting at least with \cite{kri75}, we believe that truth should be \emph{type-free}; the truth predicate (provably) applies to sentences containing itself. 

A type-free notion of truth, as it is well-known, calls for an account of the Liar and related paradoxes. We will shortly outline the solution implicit in our theories. Independently of the details, though, the view we put forward is \emph{strongly classical}: not only the external logic of the theories should be classical, but we require  the logic of their internal theories\footnote{For a theory $T$, its internal theory is defined as $\{ A \sth \text{`$\corn{A}$ is true' is provable in $T$}\}$.} to be classical as well. 
Our reasons for countenancing a classical external logic are the familiar ones. By endorsing a nonclassical theory of truth, one is bound to {sever} the links between the theory of truth and classical mathematics. This is a symptom of what McGee called `degradation of methodology' \cite[Ch.~4]{mcg91}; the standards in one's truth-theoretic theorizing shouldn't be any lower than the ones employed in theorizing about core scientific subjects.\footnote{For more on the mathematical costs of adopting a nonclassical logic, see \cite{hani18,fnd23,fie22} } 

But, we argued, the logic of the internal theory should be classical as well. The truth predicate, in natural and formal languages, is used to express important general claims. To perform its generalization role in full, the truth predicate should be \emph{fully compositional}, namely commuting with all classical logical connectives, and without any type restrictions. Full compositionality is for instance required to establish the truth of the laws of our classical external logic, such as `all instances of the law of excluded middle are true' \cite{qui90}. Compositionality is also essential to capture logical inferences involving arbitrary sentences (or propositions), what \cite{fuj22} calls {\textit{blind deductions}}. For instance, to directly formalize the argument
\begin{quote}
Everything that Gerhard says about cut-elimination is true. David asserted the negation of some of Gerhard's claims. Therefore, something David asserted is not true.
\end{quote}
one requires commutation of truth with negation in fully quantified form (as well as {and no type restriction}).

\subsection{From significance to semantic (in)sensitivity} A type-free, strongly compositional notion of truth requires a principled account of how the truth predicate can self-apply \emph{sine contradictione}. The naive principle
\beq\label{eq:truthit}
\text{for any $\vphi$, `$\vphi$ is true' is true iff $\vphi$ is true}
\eeq
is in fact inconsistent with the assumptions above. In their recent paper \cite{fuha23}, Kentaro Fujimoto and Volker Halbach propose to restrict the principle \eqref{eq:truthit} to what they call \emph{determinate} sentences, while retaining a fully compositional truth predicate. The proposed restriction of \eqref{eq:truthit}, in turn, entails the restricted $\T$-schema
\beq\label{eq:rts}
\D\corn{A}\ra (\T\corn{A}\lra A),
\eeq
{where D is a determinateness predicate.}

Several theories satisfying \eqref{eq:rts} can be found in the context of Kripke-Feferman truth. Well-known examples are Feferman's $\kf$  and $\msf{DT}$ -- from \cite{fef91} and \cite{fef08}, respectively. Such theories rely on a Russellian notion of range of significance for the interpretation of $\D$, which in the context of Feferman's theories can be further analyzed as the sentences which \emph{possess a classical truth value} relative to the intended semantics:  
    \begin{quote}
    I agree with Russell (1908) that every predicate has a domain of significance, and it makes sense to apply the predicate only to objects in that domain. In the case of truth, that domain $\D$ consists of the sentences that are meaningful and determinate, that is, have a definite truth value, true or false. $\D$ includes various but not necessarily all grammatically correct sentences that involve the notion of truth itself \cite[p.~206]{fef08}.
    \end{quote}

A similar picture is endorsed by Reinhardt, in \cite{rei86,rei85}, who labels classically true or false sentences `significant'. A well-known feature of the restriction of the $\T$-schema to determinate sentences, already extensively discussed in \cite{rei86},\footnote{See also \cite{bac15}. One applies basic logical steps to a sentence `saying of itself' that it's either not determinate or not true.} is that any such restriction would entail the existence of sentences that are provable and yet (provably) not determinate. 

Fujimoto and Halbach{, however,} reject the Russellian view and claim that truth can be (provably) applied to sentences that are not determinate. They suggest moving from the notions of range of significance to a notion of \emph{semantic sensitivity}: `some sentences, including liar sentences, are sensitive to the addition of another layer of truth. Stacking an additional layer of truth onto the liar sentence will change its semantic status; but that does not mean that the truth predicate cannot meaningfully be applied to it' \cite[p.~256]{fuha23}. The shift from `significance' to `sensitivity' has the effect of weakening the Reinhardt-Bacon issue; it is not problematic to prove sentences that are (provably) semantically sensitive. 

The shift from significance to semantic sensitivity is, in our view, a valuable step forward. Fujimoto and Halbach, however, only devote a few remarks to the notion. We sketch below our own interpretation of semantic sensitivity in the context of type-free truth, and employ it to highlight the benefits of the theories proposed in the paper. 

Another distinctive feature of Fujimoto and Halbach's account is that determinateness cannot be explained in terms of truth. Their theories feature a \emph{primitive} determinateness predicate as well as a primitive truth predicate. Fujimoto and Halbach suggest that this has to be so; in discussing Feferman's notion of determinateness from \cite{fef08}, they write `[Feferman] took $\D$ as definable in terms of truth by stipulating $\D x :\lra \T x \vee \F x$; but this definition does not yield the desired properties of $\D$ in our theory and we consequently introduce determinateness as primitive notion' \cite[p.~257]{fuha23}. 

As some of our results will indicate, the primitive nature of the determinateness  predicate significantly complicates the semantic analysis of the truth-theoretic language. Fortunately, as we will show, virtually all of the theoretical benefits sought for by Fujimoto and Halbach can be met by defining determinateness in terms of truth. 

\subsection{Logic and Semantics}

A formal theory of type-free truth plays, in our account, a double role. It provides general laws for truth that are used in several theoretical contexts, including semantic theorizing. As a result of the discussion above, these must include compositional axioms in combination with classical (external and internal) logic. In addition, the theory of truth provides a formalization of clear semantic rules to be applied to a  well-defined and sufficiently comprehensive class of sentences of the (type-free) truth-theoretic language. Because of paradox, such rules may not coincide with the classical semantic rules. We will see that this is precisely what the theories introduced below will be able to do. 

The semantic rules that we identify involve decomposition of semantic values according to the familiar Strong-Kleene truth conditions, including full disquotation for the truth predicate. In the class of models we privilege, the natural collection of sentences to which these rules unrestrictedly apply are the Kripke-determinate ones, that is sentences that are in the extension of the relevant consistent fixed-point model in the sense of \cite{kri75}.\footnote{However, we will also consider classes of models with inconsistent extensions in what follows.} Relative to each such model, these are the \emph{determinate} sentences, and the $\T$-schema holds unrestrictedly for them. Relative to the minimal fixed point, the determinate sentences are simply the grounded ones. 

As a consequence, just like in Fujimoto and Halbach's picture, the determinate sentences are semantically insensitive in the sense that their semantic status is not affected by one or more iteration of the truth predicate on the sentence. However, unlike what happens in Fujimoto and Halbach's theory, we can provide fully general and uniform semantic rules for the semantic analysis of the language with a type-free truth predicate. As the results of Section \ref{sec:intmod} will establish, such clearly defined rules cannot be found in Fujimoto and Halbach's approach. 

The logical role of truth is fulfilled, in our framework, primarily by full compositional principles. These cannot be semantically insensitive, as they are quantified sentences with instances that cannot be determinate. In fact, one of the novel features of our approach is the uniform combination of Kripkean fixed-point semantics with full compositionality. This is achieved by a move that we call \emph{classical closure of Kripkean truth}. Formally, this prescribes that the extension of the truth predicate contain the sentences that are satisfied by the relevant closed-off fixed point model. Conceptually, it demands that the truth-theorist engaged in semantic theorizing fully embrace the generalizing power afforded by the truth predicate.

It is worth noting that, as a direct consequence of full compositionality, some sentences happen to be in the extension of the classical closure of Kripkean truth only in virtue of the generalizing role of truth. This is for instance the case of classical logical truths. Take for instance the sentence $\lambda \vee\neg \lambda$ for $\lambda$ a liar sentence. The sentence is true, but any further iteration of truth on it cannot be.  

Ultimately, we are focusing on the classical models resulting from taking the classical closure of consistent fixed-points as extension of the truth predicate. Given the set up, full compositional axioms will be satisfied in such models. Moreover, since the determinate sentences can be readily defined as the ones that are recognized as true or false (i.e.~have a true negation) in the classical closure of Kripkean truth, a sentence will be determinate precisely when the sentence stating that it is true or false is itself true in the model. The model thus provides a strongly classical environment for the study of a traditional notion of determinateness. 


\subsection{Axiomatization and $\nat$-categoricity}

The paper is primarily concerned with axiomatic theories of truth capturing the intended semantics. Without axiomatization, semantic constructions are typically hard to pin down. 
Especially in the context of classical theories, a robust sense in which an axiomatic theory captures a semantic construction is given by the notion of $\nat$-categoricity introduced by \cite{fial15}. Let $\mathscr{E}\subseteq \mathscr{P}(\omega)$ be a collection of extensions of the truth predicate; a theory $T$ $\nat$-categorically axiomatizes $\mathscr{E}$ just in case
\[
\text{$(\nat,S)\vDash T$ iff $S\in \mathscr{E}$}.
\]

Another remarkable feature of the theories we introduce in the paper is that they turn out to be $\nat$-categorical with respect to the models described in the previous subsection. Specifically, a structure $(\nat,S)$ models our theories if and only if $S$ is the classical closure of a consistent Kripkean fixed point. 
The $\nat$-categoricity of our theories marks another key difference with Fujimoto and Halbach's approach afforded by the defined determinateness predicate. As the results below will indicate, while some models of Fujimoto and Halbach's theories are akin to the classical closure of Kripkean truth, some others are not. 

Therefore, there is a precise sense in which the conception of truth outlined in this section -- and embodied in the class of $\omega$-models just described -- corresponds to the laws of truth we propose and study in the paper.

\subsection{Summary of the Main Results}

From the technical side, we study theories that are closely related to the ones studied by Fujimoto and Halbach in their \cite{fuha23}, but with important difference that will be highlighted in due course. Fujimoto and Halbach mainly focus on two systems, $\cd$ and $\cdplus$ (see Dfn.~\ref{dfn:cdplus}). $\cdplus$ is obtained by strengthening the axiom
\[
\tag{T2} \forall t(\D \val{t}\ra \T\D t)
\]
of $\cd$ with the (plausible) biconditional
\[
\tag{T2$^+$} \forall t(\D \val{t}\lra \T\D t).
\]
As we will explain shortly, in the paper we focus exclusively on $\cdplus$. Here is a summary of the main results. 
\begin{enumerate}[label=(\alph*)]
\item We introduce a theory of truth, $\cdplust$, whose axioms are the axioms of $\cdplus$ except that the determinateness predicate $\D x$ is defined as $\T\T x\vee\T\F x$ (Dfn.~\ref{dfn:cdplust}). We show that the theory is consistent by displaying a natural class of $\omega$-models for it (Prop.~\ref{prop_cons-mod}), and that the theory is mutually reducible, in a strong sense, with $\cdplus$. 
\item We show that $\cdplust$ can be reaxiomatized as what we will call $\kf$'s classical closure ($\ckf$, Dfn.~\ref{dfn:ckf}), i.e.~the fully compositional theory of type-free truth (with no axioms for determinateness) obtained by {axiomatizing $\kf +\msf{CONS}$} within the scope of a truth predicate satisfying unrestricted  compositional principles (Prop.~\ref{prop:ckfcdt}). A particularly nice feature of $\ckf$ (and hence of $\cdplust$) is its $\nat$-categoricity with respect to the classical closures of consistent Kripkean fixed points (Prop.\@ \ref{pr_Ncat}). This feature is not available for $\cdplus$. 
\item While studying $\cdplust$ and $\ckf$, we provide additional results on \emph{both theories} $\cd$ and $\cdplus$: we show that they cannot prove {most} of the key axioms of $\ckf$; some of these axioms were proved to be conservative over $\cdplus$ in \cite{fuha23}. Our results show that the behaviour of $\cd$ and $\cdplus$'s truth predicate, while fully classical and compositional on the surface, becomes highly irregular within two or more layers of truth.
\item We show that our results are stable under a dual definition of $\D x$, which is suitable for  a complete (but not consistent) truth predicate. Specifically, we introduce the theory $\cdplustcp$, whose axioms are those of $\cdplus$ except that $\D x$ is defined as $\neg\T\T x\vee\neg\T\F x$, and show that it is mutually reducible with $\cdplus$. We also show that $\cdplustcp$ can be reaxiomatized as the classical closure of $\kf +\msf{COMP}$.
\end{enumerate}
{The reason for focusing on $\cdplus$ is that the natural class of models referred to in (a) and the axiomatizations mentioned in (b)-(d) license all principles of $\cdplus$.\footnote{It is unclear whether a strategy of the kind we employ to define the truth predicate of $\cdplus$ is available for $\cd$.}}


\subsection{Technical Preliminaries}
We work over the language of arithmetic $\lnat$, which extends the standard signature $\{\ovl{0},\mrm{S},+,\times\}$ with finitely many function symbols for elementary syntactic operations (see next paragraph for more details).  Let $\lt:=\lnat\cup\{\T\}$ and $\ld:=\lt\cup\{\D\}$, where $\T$ and $\D$ are unary truth- and determinateness predicates, respectively.  We fix a canonical G\"odel numbering of $\ld$-expressions and a formalization of syntactic notions and operations as it can be found, for instance, in \cite{hapu98}. Following standard practice, we take Peano arithmetic to be the first-order theory in which this formalization is carried out, although of course much weaker systems would suffice. We write $\pat$ for the theory obtained by formulating the axiom of $\pa$ in the expanded language $\lt$.  

The additional function symbols of $\lnat$ include a symbol for the standard numeral function $\mrm{num}(x)$ sending a number to the code of its numeral. A standard dot-notation to denote such symbols: 
	\begin{IEEEeqnarray*}{LL}
		\textsc{Operation} \hspace*{40mm} & \textsc{Function in $\ld$} 
		\\
		\#t, \#s\mapsto\#(t= s) & \ud=
		\\
		\#\varphi\mapsto\#(\neg\varphi) & \ng
		\\
		\#\varphi,\#\psi \mapsto\#(\varphi\land\psi) & \ud\land
		\\
		\#v_k,\#\varphi \mapsto\#(\forall v_k\varphi) & \ud\forall
		\\
		{\#t \mapsto \#\T(t)} & \ud\T\\
         {\#t \mapsto \#\D(t)} & \ud\D
	\end{IEEEeqnarray*}
 Using the conventions above, we define falsity as true negation, that is $\F x:\lra\T\ng x$. 
We also take the following $\lnat$-predicates to abbreviate the equations for the (elementary) characteristic function for such sets:
	\begin{enumerate}
        \item[] $\mrm{Term}(x) \ (\mrm{Cterm}(x)) := x$ is the G\"odel number of a (closed) term;
        \item[] $\mathrm{Fml}_{\mc L}^n(x) \ \ (\mathrm{Sent}_{\mc L}(x)) := x$ is the G\"odel number of a $\mc L$-formula with at most $n \ (0)$ free distinct variables.
    \end{enumerate}
As in Halbach's monograph \cite{hal14}, we will employ a functional notation $\val{x}$ abbreviating the formula representing in $\lnat$ the evaluation function for closed terms of $\lnat$ and is such that $\val{\corn t} = t$ for closed terms $t$.

Roman upper case letters $A, B, C,\ldots$ range over formulae of $\ld$. Greek lower case letters $\vphi,\psi,\xi,\ldots$ will be used to abbreviate quantification over formal sentences, {while $\vphi(v),\psi(v),\xi(v),\ldots$ will be used to abbreviate quantification over formal formulae with at most $v$ free}. Also, for the sake of readability, we suppress codes, dots, and Quine corners when there is no danger of confusion. Similar conventions apply to $\lt$. The expression $e(t/v_k)$, sometimes abbreviated as $e(t)$, represents the syntactic substitution of a term $t$ for a variable $v_k$ in an expression $e$. So for example, 
the expressions 
	\begin{IEEEeqnarray*}{LCL}
		(\forall\vphi,\psi\colon \ld) (\T(\vphi\land\psi) \lra\T\vphi\land\T\psi)
        \\
        \forall t (\T\T t\ra \T t^\circ)
        \\
        (\forall \vphi(v)\colon \ld) (\T(\forall v\vphi)\lra \forall t\,\T\vphi(t/v))
	\end{IEEEeqnarray*}
are short for, respectively
	\begin{IEEEeqnarray*}{LCL}
		\forall x\forall y (\mrm{Sent}_{\ld}(x)\land \mrm{Sent}_{\ld}(y) \ra (\T(x \subdot\land y) \lra\T x\land\T y))
        \\
        \forall x (\mrm{Cterm}(x)\ra \T\ud\T x\ra \T x^\circ).
        \\
        \forall x\forall y \big(\mathrm{Fml}_{\ld}^1(x)\land\mrm{Var}(y) \ra ( \T\subdot\forall yx \lra \forall z (\mrm{Cterm}(z) \ra \T x(z)))\big)
	\end{IEEEeqnarray*}
  We will also write $A\in X$ instead of $\#A\in X$, or, {for $t$ ranging over closed terms, $\T\T\T t$ instead of $\T(\mrm{\ud\T}\mrm{num}\ud\T t)$, etc}.

We use standard notions of \textit{relative translation} and \textit{relative interpretation} as it can be found, for instance, in \cite{vis06,hal14}.  We use the notions of \emph{$\lnat$-translations} and \emph{$\lnat$-interpretations}: an $\lnat$-translation is a translation between theories in $\mc{L}_1, \mc{L}\supseteq \lnat$ which does not relativize quantifiers and collapses into the identity function when restricted to $\lnat$. An $\lnat$-interpretation is an $\lnat$-translation that preserves provability in the standard way. Indeed, the special case of $\lnat$-interpretability in which $\mc{L}_1, \mc{L}\supseteq \lnat$ expand the arithmetical signature with truth predicates is the notion of truth-definability in the sense of \cite{fuj10}.

\section{$\cdplus$}

Fujimoto and Halbach introduce the theory $\cdplus$ in \cite{fuha23}. 
\begin{dfn}\label{dfn:cdplus}
$\cdplus$ consists in the following extension of $\pa$ in the language $\ld$ (where $\D$ and $\T$ are allowed to appear in induction):
\begin{align}
	\tag{$\T1$}\label{eq:t1}&\forall s\forall t(\T(s = t) \lra \val{s}=\val{t})\\
    \tag{$\T2^+$}\label{eq:t2}&\forall t (\mathrm{D}\val{t}\lra \T  \D t)\\
	\tag{$\T3$}\label{eq:t3}&\forall t(\D\val{t}\ra (\T \T t\lra \T \val{t}))\\
	\tag{$\T4$}\label{eq:t4}&(\forall \vphi\colon \ld) (\T(\neg \vphi)\lra \neg \T\vphi)\\
	\tag{$\T5$}\label{eq:t5}&(\forall \vphi,\psi)\colon \ld (\T(\vphi\land \psi)\lra \T\vphi\land \T\psi)\\
	\tag{$\T6$}\label{eq:t6}&(\forall \vphi(v)\colon \ld) (\T(\forall v\vphi)\lra \forall t\,\T\vphi(t/v))\\
	\tag{$\D1$}&\label{eq:d1}\forall s\forall t \,\D(s =t)\\
	\tag{$\D2$}\label{eq:d2}&\forall t (\D\T t\lra \D\val{t})\\
	\tag{$\D3$}\label{eq:d3}&\forall t(\D \D t\lra \D\val{t})\\
	\tag{$\D4$}&\label{eq:d4}(\forall \vphi\colon\ld)(\D(\neg \vphi)\lra \D\vphi)\\
	\tag{$\D5$}&\label{eq:d5}(\forall \vphi,\psi\colon \ld) \big(\D(\vphi\land\psi)\lra ((\D \vphi\land \D \psi)\vee (\D\vphi\land \T\neg\vphi)\vee(\D \psi\land \T\neg\psi))\big)\\
	\tag{$\D6$}\label{eq:d6}&(\forall \vphi(v)\colon\ld) (\D(\forall v\vphi)\lra (\forall t \,\D\vphi(t/v)\vee \exists t \,\D\vphi(t/v)\land \T\neg\vphi(t/v)).\\
    \tag{$\mrm{R}1$}\label{eq:r1}&(\forall \vphi(v)\colon \ld)\forall s\forall t\big(\val s= \val t \ra (\T\vphi(s) \lra \T\vphi(t))\big)\\
    \tag{$\mrm{R}2$}\label{eq:r1}&(\forall \vphi(v)\colon \ld)\forall s\forall t\big(\val s= \val t \ra (\D\vphi(s) \lra \D\vphi(t))\big)
\end{align}
\end{dfn}

Fujimoto and Halbach establish the consistency of $\cdplus$ -- and, of course, of $\cd$ -- by exhibiting an $\omega$-model. Since the construction will be relevant later on, we repeat it here. Let $\mathcal{D}(x)$ be the positive inductive definition associated with the (right-to-left direction) of the $\D$-axioms $\D1$-$\D6$ \cite[\S4]{fuha23}. Let 
\[
	\Gamma_Y(X)=\{ n \sth (\nat,X,Y)\vDash \mc{D}(\ovl{n})\},
\]
and define \label{semcdp}
\begin{align*}
	&D_0:=\varnothing&& T_0:=\varnothing\\
	&D_{\alpha+1}:=\Gamma_{T_\alpha}(D_\alpha)&& T_{\alpha+1}:=\{A \in \ld\sth (\nat,D_\alpha,T_\alpha)\vDash A\}\\
	&D_\lambda:=\bigcup_{\beta<\lambda} D_\beta &&T_\lambda :=\bigcup_{\beta<\lambda} \big( D_\beta\cap T_\beta\big)
\end{align*}

Then $D_{\omega_1}$ is a fixed point of $\Gamma_{T_{\omega_1}}$. Let us denote $D_{\omega_1}$ and $T_{\omega_1}$ by $D_\infty$ and $T_\infty$, respectively. Let:
\[
	\mbb T_\infty=\{ A\in \ld \sth (\nat, D_\infty,T_\infty)\vDash A\}. 
\]

\begin{lemma}[{\cite[Theorem 4.8]{fuha23}}]
$(\nat, D_\infty,\mbb T_\infty)\vDash \cdplus$. 
\end{lemma}

The proof-theoretic analysis of $\cdplus$ is quite straightforward. It reduces the system to a theory resulting from a combination of typed truth with standard Kripke-Feferman systems. It employs $\ctkfc$, that is the result of enriching the Kripke-Feferman theory with consistency with a \emph{typed, Tarskian} truth predicate $\mbf{T}$ on top, but a version of $\kf$ with the completeness axiom would work as well.
\begin{dfn}\label{df:kf}
$\kf$ is the system in $\lt$ extending $\pat$ with the following axioms:
\begin{align*}
\tag{$\kf 1$}& \forall s\forall t(\T(s=t)\lra \val{s}=\val{t})\land \forall s\forall t(\T(s\neq t)\lra \val{s}\neq \val{t})\\
\tag{$\kf 2$}& \forall t (\T\T t\lra \T\val t)\land \forall t (\F\T t\lra \F\val{t}\vee\neg \mrm{Sent}(\val{t}))\\
\tag{$\kf 3$} &(\forall \vphi\colon \lt)(\F \neg \vphi\lra \T\vphi)\\
\tag{$\kf4$}&(\forall \vphi,\psi\colon \lt)((\T(\vphi\land \psi)\lra \T\vphi\land \T\psi)\land (\F(\vphi\land \psi)\lra \F\vphi \vee \F\psi))\\
\tag{$\kf5$}&(\forall \vphi(v)\colon \lt)\big((\T(\forall v\vphi)\lra \forall t \,\T\vphi(t/v))\land (\T(\exists v\vphi)\lra \exists t\,\T\vphi(t/v))\big)
\end{align*}
$\kfcon$ is obtained by adding to $\kf$ the axiom
\[
\tag{$\msf{Cons}$} (\forall \vphi\colon \lt)(\F \vphi \ra \neg \T\vphi);
\]
$\kfcom$ is obtained by adding to $\kf$ the converse claim
\[
\tag{$\msf{Comp}$} (\forall \vphi \colon \lt)(\neg\T\vphi \ra \F \vphi).
\]
\end{dfn}

\begin{dfn}[$\ctkfc$]\label{df:dfctkf}
We expand $\lt$ with an additional truth predicate $\mbf{T}$ and call the resulting language $\lctkc$. $\ctkfc$ is the theory in $\lctkc$ extending $\kf+\msf{CONS}$ with the following axioms:
\begin{align}
\tag{$\mbf{T} 1$}& (\mbf{T}(s=t)\lra \val{s}=\val{t})\land \mbf{T}\T t\lra  \T \val{t}\\
\tag{$\mbf T 2$}& (\forall \vphi\colon \lt)(\mbf T(\neg \vphi)\lra \neg \mbf T\vphi)\\
\tag{$\mbf T 3$}& (\forall \vphi,\psi \colon \lt)(\mbf T(\vphi\land \psi)\lra \mbf T\vphi \land \mbf T\psi)\\
\tag{$\mbf T 4$}&(\forall \vphi(v)\colon \lt)(\mbf{T}(\forall v\vphi)\lra \forall t \,\mbf T(\vphi(t/v)))
\end{align}
\end{dfn}

We assume a standard notation for ordinals below the Feferman-Sch\"utte ordinal $\Gamma_0$ (see e.g.~\cite[ch.~2]{poh09}). In particular, we denote with $\varepsilon_\alpha$ the (code of) the $\alpha^{th}$ fixed point of the function $\lambda x.\,\omega^x$. 

\begin{lemma}[{\cite[Theorem 7.12]{fuha23}}] $\cdplus$ is mutually $\lnat$-interpretable with $\ctkfc$, and therefore with ramified truth/analysis up to any $\alpha<\varepsilon_{\varepsilon_0}$. 
\end{lemma}

\begin{proof}[Proof Idea]
To interpret $\cdplus$ in $\ctkfc$, external occurrences of $\T$ are translated by the Tarskian $\mbf{T}$, whereas internal occurrences of $\T$ are translated homophonically as $\T$. $\D$ becomes $\T \cup \mrm{F}$ (the sentences with `a classical truth-value' in $\kfcon$). 
Conversely, to interpret $\ctkfc$ in $\cdplus$, the Tarskian $\mbf{T}$ is translated as $\T$, and the $\kfcon$ truth predicate $\T$ is instead translated as the `determinate truths' of $\cdplus$, i.e.~$\T \cap \D$. 

\end{proof}

\section{$\cdplust$ and its Intended Models}\label{sec:intmod}

Our first result consists in displaying a natural class of models for type-free truth, the \emph{classical closures of (consistent) Kripkean fixed points}. We will also show that such structures model {the axioms of} Fujimoto and Halbach's $\cdplus$. Crucially, in these models, $\D$ is defined in terms of the notion of Kripkean determinateness. In the next section, we will provide direct axiomatizations of such structures in the language $\lt$.   

We recall the arithmetical operator associated with the Strong-Kleene version of the fixed-point semantics. 
\begin{dfn}[$\mrm{K}$-jump]
For any $X\subseteq \mrm{Sent}_{\lt}$, 
\begin{IEEEeqnarray*}{rCl}
 n\in \mathscr{K}(X)& :\lra & n\in \mrm{Sent}_{\lt}\land 
  \\
  && \big(\exists s\exists t (n=(s=t)\land \val{s} = \val{t}) \vee
  \\
  && \big(\exists s\exists t (n=(s\neq t)\land \val{s} \neq \val{t}) \vee\\
  && \big(\exists t (n=(\T t)\land \val{t}\in X) \vee\\
  && \exists t (n=(\neg \T t)\land \neg\val{t}\in X\vee \val{t}\in \omega\setminus \mrm{Sent}_{\lt})\vee\\
  && \exists \vphi(n=\neg\neg \vphi \land \vphi \in X)\vee\\
  &&\exists \vphi,\psi(n=(\vphi\land \psi)\land \vphi \in X\land \psi\in X)\vee\\
  && \exists \vphi,\psi(n=\neg(\vphi\land \psi)\land \neg \vphi \in X \vee \neg \psi\in X)\\
  && \exists \vphi(v)(n=(\forall v\vphi)\land \forall t \,\vphi(t/v)\in X)\\
  && \exists \vphi(v)(n=\neg \forall v\vphi\land \exists t \,\neg\vphi(t/v)\in X)\big)
\end{IEEEeqnarray*}
\end{dfn}
\noindent $X$ occurs positively in $\mathscr{K}(X)$ and the $\mrm{K}$-jump is monotone. Therefore, it will have fixed points. Note that fixed-points of $\msc{K}$ are closed under Strong Kleene logic and are such that $(\neg)\vphi\in X$ iff $(\neg)\T\vphi\in X$.\footnote{For more details on fixed-point semantics, we refer to \cite[\S15.1]{hal14} and \cite[\S\S4-5]{mcg91}.}

In this section and the next we are interested in \emph{consistent fixed-points $\msc{K}$} (and theories that are sound with respect to them), that is sets $X$ such that $X=\msc{K}(X)$ and such that there's no $\vphi$ such that $\vphi\land \neg \vphi \in X$; we will consider other fixed-points in later section. The reason for this is both technical and conceptual. Technically, consistency delivers a simpler definition of determinateness and provides a basis for results about complete models and theories presented in later sections. Conceptually, we believe that a consistent extension of truth provides a more attractive notion of determinateness, as well as being in continuity with the above-mentioned works on truth-theoretic determinateness by Kripke, Reinhardt and Feferman.

For a consistent fixed-point $X$, let
\beq
    \mbb{T}_X=\{A\in \lt \sth (\nat,X)\vDash A\}.
\eeq
Incidentally, this set is considered by Fujimoto and Halbach as a model of $\ctkfc$ \cite[p.~251]{fuha23}. We will now show that this set can be used to provide a direct model construction for the principles of $\cdplus$, and not only for a {typed} theory of truth interpreting $\cdplus$.

\begin{dfn}[$\cdplust$]\label{dfn:cdplust}
$\cdplust$ is the theory in $\lt$ whose axioms are the axioms of $\cdplus$, but where $\D$ is now \emph{defined} in terms of $\T$ as
{\begin{IEEEeqnarray}{L}
\D x:\lra \T\subdot\T \mrm{num}(x)\lor \T\subdot\F\mrm{num}(x)\;\;\;\text{(abbr.~$\T\T x\vee \T\F x$)}. \label{eq_Ddf-cons}
\end{IEEEeqnarray}}
\end{dfn}

$\cdplust$, just like $\cdplus$, delivers the intended restriction to the $\T$-schema to determinate sentences. 
\begin{obse}
For any $A\in \lt$: $\cdplust\vdash \D\corn{A}\ra (\T\corn{A}\lra A)$. 
\end{obse}
\begin{proof}
By external induction on the complexity of $A$. Note that the case in which $A$ is $\T t$ is an axiom of $\cdplust$. 
\end{proof}

Our next goal is to show that, given a consistent fixed-point $X$, the structure $(\nat,\mbb T_X)$ is a model of $\cdplust$. Since $(\nat,X)\vDash\kfcon$, and since $\kfcon\vdash\T\corn A\ra A$, we have

\begin{fact}\label{fact_auxiliary}
	For any consistent fixed-point $X=\msc K(X)$, $(\nat,X)\vDash\T\corn A\ra A$.
\end{fact}

\begin{prop}\label{prop_cons-mod}
	$(\nat,\mbb T_X)\vDash \cdplust$. 
\end{prop}
\begin{proof}
By induction on the length of the proof in $\cdplust$. We verify some key axioms, noting that T1, T4-T6 are immediate by definition.

\medskip
\noindent
T2$^+$: $(\nat,\mbb{T}_{X})\vDash\T\T\vphi \lor \T\F\vphi$ iff $(\nat,X)\vDash\T\vphi \lor \F\vphi$ iff $\{\vphi,\neg\vphi\}\cap X\neq\vno$ iff $\{\T\vphi,\F\vphi\}\cap X\neq\vno$ (by the fixed-point property), iff $(\nat,X)\vDash\T\T\vphi\lor\T\F\vphi$ iff $(\nat,\mbb{T}_X)\vDash\T(\T\T\vphi\lor\T\F\vphi)$.

\medskip
\noindent
T3: Assume $(\nat,\mbb{T}_{X})\vDash\T\T\vphi \lor \T\F\vphi$, which is the case iff $\{\vphi,\neg\vphi\}\cap X\neq\vno$. To show $(\nat, \mathbb{T}_X)\vDash\T\T\vphi\ra\T\vphi$, we reason as follows, letting $\vphi=\corn A$:\footnote{
	We assume that $\vphi$ denotes a sentence without loss of generality, as $X\subseteq\mrm{Sent}_{\mathcal{L}_\T}$.
}
	\begin{IEEEeqnarray*}{LL}
		(\nat, \mathbb{T}_{X})\vDash\T\T\vphi \hspace{3mm} & \text{ iff }
		\\
		(\nat,{X})\vDash\T\vphi & \text{ \textit{hence}, by Fact \ref{fact_auxiliary}},
        \\
        (\nat, X)\vDash A, & \text{ iff}
        \\
        (\nat, \mathbb{T}_{X})\vDash\T\vphi.
	\end{IEEEeqnarray*}
 We notice that -- due to the fact that we are reasoning in a consistent fixed-point model -- the assumption $\D\vphi$ has not been employed in this part of the argument. 

Conversely, to show that $(\nat, \mathbb{T}_X)\vDash\T\vphi\ra\T\T\vphi$, assume $(\nat, \mathbb{T}_{X})\vDash\T\vphi$, {which is equivalent to $(\nat, X)\vDash A$}. Towards a contradiction, suppose $(\nat, \mbb T_X)\not\vDash\T\T\vphi$, hence $(\nat, X)\not\vDash\T\vphi$. Then, since $\{\vphi,\neg\vphi\}\cap X\neq\vno$ by assumption,  $(\nat, X)\vDash\F\vphi$, and therefore $(\nat, X)\vDash \neg A$ by Fact~\ref{fact_auxiliary}, which contradicts $(\nat, X)\vDash A$.

\medskip
\noindent
D3: We need to show that $(\mbb{N}, \mbb T_X)\vDash\T\T\vphi\lor\T\F\vphi$ is equivalent to
	\begin{IEEEeqnarray*}{RCL}
		(\mbb{N}, \mbb T_X)\vDash\T\T \big(\T\T\vphi\lor\T\F\vphi\big) \lor 
		\T\F \big(\T\T\vphi\lor\T\F\vphi\big).
	\end{IEEEeqnarray*}
The latter is equivalent to
	\begin{IEEEeqnarray*}{RCL}
		\underbrace{(\mbb{N}, X)\vDash\T\big(\T\T\vphi\lor\T\F\vphi\big)}_{d_1}
		\ &\text{ or }& \ 
		\underbrace{(\mbb{N}, X)\vDash\F\big(\T\T\vphi\lor\T\F\vphi\big)}_{d_2}.
	\end{IEEEeqnarray*}

Since $d_1$ and $d_2$ are equivalent to, respectively, $(\mbb{N}, X)\vDash\T\vphi\lor\F\vphi$ and $(\mbb{N}, X)\vDash\F\vphi\land\T\vphi$, by consistency of $X$ it can be observed that their disjunction is equivalent to $d_1$, which is equivalent to $(\mbb{N}, \mbb T_X)\vDash\T\T\vphi\lor\T\F\vphi$, as required.

\medskip
\noindent
D5: Assume $(\mbb N,\mbb T_X)\vDash\T\T(\vphi\land\psi)\lor\T\F(\vphi\land\psi)$. If $(\mbb N,\mbb T_X)\vDash\T\T(\vphi\land\psi)$, then $(\mbb N,\mbb T_X)\vDash\T\T\vphi\land\T\T\psi$ by the closure properties of $X$, hence both $\vphi$ and $\psi$ are determinate. If $(\mbb N,\mbb T_X)\vDash\T\F(\vphi\land\psi)$ then $(\mbb N,\mbb T_X)\vDash\T\F\vphi\lor\T\F\psi$ by the properties of $X$. If $(\mbb N,\mbb T_X)\vDash\T\F\vphi$, then $(\mbb N, X)\vDash\F\vphi$ and hence $(\mbb N, X)\vDash\neg A$ by Fact~\ref{fact_auxiliary}, for $\vphi=\corn A$. This in turn yields $(\mbb N,\mbb T_X)\vDash\T\F\vphi\land\F\vphi$, hence $\vphi$ is determinate and false. Similarly if $(\mbb N,\mbb T_X)\vDash\T\F\psi$.

Conversely, if $(\mbb N,\mbb T_X)\vDash(\T\T\vphi\lor\T\F\vphi)\land(\T\T\psi\lor\T\F\psi)$, then we get immediately $(\mbb N,\mbb T_X)\vDash\T\T(\vphi\land\psi)\lor\T\F(\vphi\land\psi)$ by closure properties of $X$. If $(\mbb N,\mbb T_X)\vDash(\T\T\vphi\lor\T\F\vphi)\land\F\vphi$, then $(\mbb N,X)\vDash(\T\vphi\lor\F\vphi)\land\neg A$, which by consistency of $X$ is equivalent to $(\mbb N,X)\vDash\F\vphi$, hence $(\mbb N,X)\vDash\F(\vphi\land\psi)$ for any $\psi$. It follows that $(\mbb N,\mbb T_X)\vDash\T\F(\vphi\land\psi)$. Similarly if we assume $(\mbb N,\mbb T_X)\vDash(\T\T\psi\lor\T\F\psi)\land\F\psi$.
\end{proof}

It is clear that $\cdplust$ can interpret $\cdplus$. 
\begin{obse}
There is an $\lnat$-interpretation of $\cdplus$ in $\cdplust$. 
\end{obse}
\begin{proof}
We employ the recursion theorem to define an $\lnat$-translation $\delta \colon \ld \to \lt$ which systematically replaces $\D$ with its definition from \eqref{eq_Ddf-cons}.
In more detail, the recursion theorem for primitive recursive functions\footnote{See for instance \cite[\S11.2]{rog67}.} can be employed to define an $\lnat$-translation $\delta \colon \ld \to \lt$ such that
\[
    (\D x)^\delta :\lra \T\subdot \T  \mrm{num}(\subdot\delta (x))\vee \T\subdot\F \mrm{num}(\subdot\delta (x)).
\]
The expression $\subdot\delta(x)$ abbreviates the formula representing  $\delta$ in $\lnat$. 
The verification that 
\beq
\cdplus \vdash A \text{ only if } \cdplust\vdash A^\delta
\eeq
is immediate given some basic syntactic facts, provable in a subtheory of $\pa$, including:
\begin{align}
& \forall x (\mrm{Sent}_{\ld}(x)\ra \mrm{Sent}_{\lt}(\subdot\delta(x))). 
\end{align}
\end{proof}

It will follow from the identity of $\ckf$ and $\cdplust$ (Prop.~\ref{prop:ckfcdt}) that  $\cdplus$ and $\cdplust$ are mutually $\lnat$-interpretable. {It would be too hasty to think, however, that $\cdplus$ and $\cdplust$ are ``notational variants''}; $\cdplus$ cannot define $\D$ in the manner prescribed by  $\cdplust$. This can be seen from the following observation. With reference to Fujimoto and Halbach's semantic construction for $\cdplus$ described on page \pageref{semcdp}, let 
\[
    T_B:=T_\infty\cup \{B\},
\]
for $B$ an $\ld$-sentence, and 
\[
    \mbb{T}_B:=\{ A \in \ld \sth (\nat,D_\infty,T_B)\vDash A\}. 
\]
\begin{lemma}\label{lem:cmdefd}
If $B \notin D_\infty$, then $(\nat, D_\infty,\mbb{T}_B)\vDash \cdplus$. 
\end{lemma}
\begin{proof}
	We first show, by induction on $A$, that if $A\in D_\infty$, then
	\begin{IEEEeqnarray}{R}
		(\mbb N, D_{\infty}, T_{\infty})\vDash A
		\, \text{ iff } \,
		(\mbb N, D_{\infty}, T_B)\vDash A. \label{eq_aux-gamma}
	\end{IEEEeqnarray}
The crucial case is when $A=\T\vphi$: $(\mbb N, D_{\infty}, T_{\infty})\vDash\T\vphi$ iff $\vphi\in T_{\infty}$, \textit{hence} $(\mbb N, D_{\infty}, T_B)\vDash\T\vphi$. For the converse direction, $(\mbb N, D_{\infty}, T_B)\vDash\T\vphi$ iff either $\vphi\in T_{\infty}$, in which case we are done, or $\vphi=B$. However, since $\T\vphi\in D_\infty$, it follows that $\vphi\in D_\infty$, hence $\vphi\neq B$. 

Having shown \eqref{eq_aux-gamma}, we proceed with the main claim by induction on the length of proofs of $\cdplus$. We check two cases where the use of \eqref{eq_aux-gamma} is relevant, letting $\vphi=\corn A$.

\begin{itemize}
\item[-] For T3, assume $A\in D_\infty$. Then $(\mbb N, D_{\infty}, T_{\infty})\vDash A$ iff $A\in T_{\infty}$ since $(\mbb N, D_{\infty}, \mbb T_{\infty})$ models $\cdplus$. Then $(\mbb N, D_{\infty}, \mbb T_B)\vDash\T\vphi$ iff $(\mbb N, D_{\infty}, T_B)\vDash A$ iff, by \eqref{eq_aux-gamma}, $(\mbb N, D_{\infty}, T_{\infty})\vDash A$ iff $A\in T_{\infty}$ iff, since $A\neq B$, $(\mbb N, D_{\infty}, T_B)\vDash\T\vphi$ iff $(\mbb N, D_{\infty}, \mbb T_B)\vDash\T\T\vphi$.

\item[-]
For D5, the left-to-right direction is immediate. As for right-to-left direction, it suffices to consider the case where, e.g., $(\mbb N, D_{\infty}, \mbb T_B)\vDash\D\vphi\land\F\vphi$. This is the case iff $(\mbb N, D_{\infty}, T_B)\vDash\D\vphi\land\neg A$ iff, by \eqref{eq_aux-gamma}, $(\mbb N, D_{\infty}, T_{\infty})\vDash\D\vphi\land\neg A$ iff $(\mbb N, D_{\infty}, \mbb T_{\infty})\vDash\D\vphi\land\F\vphi$, hence $\vphi\land\psi\in D_\infty$ for any $\psi$, hence $(\mbb N, D_{\infty}, \mbb T_B)\vDash\D(\vphi\land\psi)$. \qedhere
\end{itemize}
\end{proof}

While $\cdplus\vdash(\forall\vphi\colon\ld)(\D\vphi \ra \T\T\vphi\lor\T\F\vphi)$, Lemma \ref{lem:cmdefd} yields that the converse implication is not provable:
\begin{corollary}\label{cor:cdpdt}
	$\cdplus\not\vdash (\forall \vphi\colon \ld)(\T\T\vphi\lor\T\F\vphi\ra\D\vphi)$.
\end{corollary}
\begin{corollary}\label{cor:cdpit}
$\cdplus$ \emph{does not} prove any of the following sentences:
\begin{enumeratei}\setlength\itemsep{1ex}
\item $\forall t( \T\T \T t \ra \T\T t)$;
\item $ \forall t( \T\T t \ra \T\T\T t)$;
\item $(\forall \vphi \colon \lt)(\T\T\vphi\lra \T\T\neg \neg \vphi)$;
\item $(\forall \vphi,\psi \colon \lt) (\T\T(\vphi\land \psi)\lra \T\T\vphi\land \T\T\psi)$;
\item $(\forall \vphi,\psi \colon \lt) (\T\F(\vphi\land \psi)\lra \T\F\vphi\lor \T\T\psi)$;
\item  $(\forall \vphi(v)\colon \lt) (\T\T\forall \vphi(v) \lra\T\forall t\T\vphi(t)$;
\item  $(\forall \vphi(v)\colon \lt) (\T\F\forall \vphi(v) \lra\T\exists t\F\vphi(t)$.
\end{enumeratei}
\end{corollary}
\begin{proof}
For (i), to construct the required countermodel use $\neg \lambda$, for $\lambda$ a liar sentence, in place of $B$ in Lemma \ref{lem:cmdefd}; the model invalidates $\T\T\T\lambda \ra \T\T\lambda$. For (ii), use $\lambda$ in place of $B$ so to invalidate $\T\T\lambda \ra \T\T\T\lambda$. The other cases can be dealt with by emploing suitable truth-teller sentences. For instance, for (iv), define parametrized truth-tellers $\tau_0$ and $\tau_1$, and replace $\tau_0\land \tau_1$ for $B$ in Lemma \ref{lem:cmdefd}. 
\end{proof}

\section{$\kf$'s Classical Closure}

The definition of $\cdplust$ is parasitic on Fujimoto and Halbach's $\cdplus$. As such, they are axiomatizations of \emph{both truth and determinateness}. 
It is then natural to ask whether there is an \emph{axiomatization of truth} that can be directly inspired by the models $(\nat, \mbb{T}_X)$ introduced above and yet deliver the required principles for determinateness. In such models, full disquotation is allowed only under two or more layers of truth (or falsity): outside those layers, fully classical principles are licensed at the expense of disquotation. {We will now answer the question positively}. We will then explain in the final section how the theories capture the conception of truth (and determinateness) outlined in the opening section.

We introduce the theory $\ckf$ (standing for the classical closure of $\kfcon$). $\ckf$ combines full compositionality ($\T1$, $\T4$-$\T6$) with Kripkean truth conditions in the inner layers of truth, {along with a disquotation principle for truth ascriptions that will ensure the consistency of the inner truth predicate}.

\begin{dfn}[$\ckf$]\label{dfn:ckf}
The axioms of the $\lt$-theory $\ckf$ are the axioms of $\pat$ together with:
\begin{align}
\tag{$\T1$} &\forall s\forall t(\T(s=t)\lra \val{s}=\val{t})\\
\tag{$\T4$} & (\forall \vphi \colon \lt) (\T(\neg \vphi)\lra \neg \T\vphi)\\
\tag{$\T5$} & (\forall \vphi,\psi\colon \lt)(\T(\vphi\land \psi)\lra \T\vphi\land \T\psi)\\
\tag{$\T6$}& (\forall \vphi(v)\colon \lt)(\T(\forall v\vphi)\lra \forall t\,\T\vphi(t/v))\\
\tag{$\T7$} &\forall s\forall t(\T\T s=t \lra \val{s}= \val{t})\\
\tag{$\T8$} &\forall s\forall t(\T\F s=t\lra \val{s}\neq \val{t})\\
\tag{$\T9$} &(\forall \vphi \colon \lt)(\T\T\vphi\lra \T\T\neg \neg \vphi)\\
\tag{$\T10$} & (\forall \vphi,\psi \colon \lt)(\T\T(\vphi\land \psi)\lra \T\T\vphi\land \T\T\psi)\\
\tag{$\T11$} & (\forall \vphi,\psi \colon \lt) (\T\F(\vphi\land \psi)\lra \T\F\vphi\lor \T\T\psi)\\
\tag{$\T12$} & (\forall \vphi(v)\colon \lt) (\T\T\forall \vphi(v) \lra\T\forall t\T\vphi(t)\\
\tag{$\T13$} & (\forall \vphi(v)\colon \lt )(\T\F\forall \vphi(v) \lra\T\exists t\F\vphi(t)\\
\tag{$\T14$} &\forall t(\T\T\T t\lra \T\T t) \\
\tag{$\T15$} &\forall t(\T\F\T t\lra \T\F t)\\
\tag{$\T$Del} &\forall t(\T\T t\ra\T\val t)\\
\tag{$\mrm{R}1$}\label{eq:r1}&(\forall \vphi(v)\colon \lt)\forall s\forall t\big(\val s= \val t \ra (\T\vphi(s) \lra \T\vphi(t))\big)
\end{align}
\end{dfn}

\begin{remark}
In the definition of $\ckf$, the version of $\kf$ that `lives' inside one layer of truth is slightly different from the one presented in Definition \ref{df:kf}. Specifically, axiom $\T15$ of $\ckf$ does not feature, internally, the extra disjunct of $\kf2$. This reformulation, albeit inessential either conceptually or from the perspective of proof-theoretic strength, is required by our arguments to obtain the identity of $\ckf$ and $\cdplust$ established in Proposition \ref{prop:ckfcdt}.
\end{remark}

We list some basic theorems of $\ckf$ that will be explicitly used below:\footnote{
	The statement \textup{\ref{eq:inta}} amounts to the consistency axiom {\sf Cons} under one additional layer of truth. Similarly to (TDel), it ensures the consistency of the inner truth predicate. We will come back to this below, \S\ref{sec:AltAx}.
 }

\begin{obse}\label{obs_ckf-theorems}
The following are derivable in $\ckf$:
	\begin{IEEEeqnarray}{L}
		\forall t(\T\F t \lra\T\T\F t) \label{eq:obs_equivalence}
        \\
        \T\F\vphi \ra \neg \T\T\vphi \ztag{TCons} \label{eq:inta}
        \\
         (\forall \vphi,\psi\colon \lt)(\T\T(\vphi\lor\psi) \lra\T\T\vphi\lor\T\T\psi) \label{eq:obs_Tdisj}
        \\
        (\forall \vphi,\psi\colon \lt)( \T\F(\vphi\lor\psi) \lra\T\F\vphi\land\T\F\psi) \label{eq:obs_Fdisj}
	\end{IEEEeqnarray}
\end{obse}

As mentioned, the theory $\ckf$ is directly inspired by the models $(\nat, \mbb T_X)$ introduced in the previous section. Our next results provide a full characterization of standard models of $\ckf$ in these terms: we show that they are exactly the classical closures of consistent fixed-points.

\begin{lemma}\label{le_NCat1}
For any consistent fixed-point $X$, $(\nat,\mbb{T}_X)\vDash \ckf$. 
\end{lemma}
\begin{proof}
    Axioms T1 and T4-T6 are immediate by definition. Axioms T7-T15 and R1 readily follows from properties of $X$. As for T16, we reason as in the argument for case T3 in the proof of Proposition~\ref{prop_cons-mod}.
\end{proof}

For any set $S$, let $S^\mathscr{I} :=\{A\sth(\mbb N, S)\vDash\T\corn{\T \corn A}\}$.
Then

\begin{lemma}\label{le_NCat2}
	Let $(\mbb N, S)\vDash\ckf$. Then $S^\mathscr{I}$ is a consistent fixed-point.
\end{lemma}

\begin{proof}
	Let $(\mbb N, S)\vDash\ckf$. By induction on the positive complexity of $A$, one shows that $A\in S^\mathscr{I}$ iff $A\in\mathscr K(S^\mathscr{I})$. For example, if $A$ is of the form $\neg\T t$, one uses T15, or if $A$ is of the form $\forall xB$, one uses T6 and T12.
\end{proof}

\begin{lemma}\label{le_NCat3}
	Let $(\mbb N, S)\vDash\ckf$. Then $S=\mbb T_{S^\mathscr{I}}$.
\end{lemma}

\begin{proof}
	Let $(\mbb N, S)\vDash\ckf$. Again by positive induction on $A$, one shows that $(\mbb N, S)\vDash\T\corn A$ iff $(\mbb N, S^\mathscr{I})\vDash A$. For example, if $A$ is $\T\corn B$, then $(\mbb N, S)\vDash\T\corn{\T\corn B}$ iff $B\in S^\mathscr{I}$ iff $(\mbb N, S^\mathscr{I})\vDash A$.
\end{proof}

From lemmata \ref{le_NCat1}, \ref{le_NCat2}, \ref{le_NCat3}, we obtain the following characterization theorem:

\begin{prop}[$\mbb N$-Categoricity]\label{pr_Ncat}
	$(\mbb N, S)\vDash\ckf$ iff $S=\mbb T_{S^\mathscr{I}}$.
\end{prop}

There's also a precise sense in which the logical strength of $\ckf$ coincides with the one of $\ctkfc$, and therefore of $\cdplus$. 
\begin{prop}\label{pr_ckf-ctcon-cdplus}
$\ckf$, $\msf{CT[\kfcon]}$, and $\cdplus$ are mutually $\lnat$-interpretable. 
\end{prop}
%
\begin{proof}
To $\lnat$-interpret $\ckf$ in $\msf{CT[\kfcon]}$, one defines an $\lnat$-translation $\tau \colon \lt \to \lctkc$ that replaces outer occurrences of $\T$ with $\mbf{T}$:
\begin{align*}
    &\tau(s=t):\lra  s=t\\
    &\tau(\T x):\lra \mbf{T}x\\
    &\tau(\neg A):\lra \neg \tau(A)\\
    & \tau(A\land B):\lra \tau(A)\land \tau(B)\\
    &\tau(\forall v A):\lra \forall x \tau (A(x/v))
\end{align*}
The translation $\tau$ is only external, and does not require applications of the primitive recursion theorem in its definition. The verification that $\tau$ is an $\lnat$-interpretation is fairly straightforward. The compositional principles $(\T1,\T4\text{-}\T6)$ follow from the definition of the translation $\tau$ and the compositional axioms for $\mbf{T}$ of $\ctkfc$. Axioms $\T7\text{-}\T13$ follow from the definition of $\tau$, $\mbf{T}1$, and the corresponding axiom of $\kf$. For instance, for $\T7$ (in non abbreviated form):
\begin{IEEEeqnarray*}{rCl+r}
\tau (\T\subdot \T \mrm{num}(s\subdot = t)\lra \val{s}=\val{t}) & \text{ iff }&\mbf{T}\subdot \T\mrm{num}(s\subdot = t)\lra \val{s}=\val{t} &  \text{def.~of $\tau$}\\
    & \text{ iff }& \T (s\subdot = t)\lra \val{s}=\val{t} & \text{by $\mbf{T}1$}
\end{IEEEeqnarray*}
The last line is an axiom of $\ctkfc$. Axiom $\mrm{R1}$ is proved by formal induction on the complexity of $\vphi(v)$. 

For $\T\mrm{Del}$, one notices that, by the definition of $\tau$ and $\mbf{T}1$, the translation of $\T\mrm{Del}$ becomes
\beq\label{eq:tra1}
\T \val{t}\ra \mbf{T}\val{t}. 
\eeq
Now either $\neg \mrm{Sent}_{\lt}(\val{t})$ or $\mrm{Sent}_{\lt}(\val{t})$.
In the latter case we can then prove \eqref{eq:tra1} by formal induction on the positive complexity of $\val{t}$. Notice that in the case in which $\val{t}$ is of form $\neg \T s$, we employ $\msf{Cons}$. If $\neg \mrm{Sent}_{\lt}(\val{t})$, then $\F\T t$ by $\kf2$, and $\neg \T\T t$ by $\msf{Cons}$, thus $\neg \T\val{t}$ and $\T\val{t}\ra \mbf{T}\val{t}$ follows trivially. 

To interpret  $\ctkfc$ in $\ckf$, we first define the internal translation $\iota \colon \lctkc\to \lt$ operating on codes of $\lctkc$-formulae as follows:
\begin{IEEEeqnarray*}{rCll}
    \iota(\vphi)&:=& \vphi, & \text{ for $A$ atomic of $\lnat$}\\
    \iota(\psi)&:=& \corn{0=1},  &\text{ for $\psi\in\mrm{Sent}_{\lctkc}\setminus\mrm{Sent}_{\lt}$} \\
    \iota(\T x)&:= & \T\subdot \iota(x)\\
    \iota(\neg \vphi)&:= & \neg \iota(\vphi)\\
    \iota(\vphi\land \psi)&:= & \iota(\vphi)\land \iota(\psi)\\
    \iota(\forall v\vphi)&:= & \forall x\iota(\vphi(x/v)).
\end{IEEEeqnarray*}
We then define the full translation $\sigma \colon\lctkc\to \lt$, which replaces $\mbf{T}$ with one single layer of $\ckf$-truth and $\T$ by two, and behaves internally according to $\iota$:
\begin{IEEEeqnarray*}{rCll}
    \sigma(A)&:\lra& A, & \text{ for $A$ atomic of $\lnat$}\\
    \sigma(\mbf{T} x)&:\lra & \T\subdot \iota(x)\\
    \sigma(\T x)&:\lra& \T\subdot \T \mrm{num}(\subdot\iota (x))\\
    \sigma(\neg A)&:\lra & \neg \sigma(A)\\
    \sigma(A\land B)&:\lra & \sigma(A)\land \sigma(B)\\
    \sigma(\forall vA)&:\lra & \forall x\sigma(A(x/v)).
\end{IEEEeqnarray*}
It remains to verify that $\sigma$ is the required $\lnat$-interpretation of $\ctkfc$ in $\ckf$.  We consider the crucial case of $\kf2$. 
 The translation of its first conjunct is 
 \beq\label{eq:demk2a}
 \T\subdot \T\mrm{num}(\subdot\T\subdot \iota( t))\lra \T\subdot \T(\subdot \iota (t))
 \eeq
 which follows immediately from $\T14$. The translation of the second conjunct of $\kf2$ is 
 \beq\label{eq:demk2b}
\T\subdot \T\mrm{num}\subdot \neg \subdot\T\subdot \iota(t)\lra \T\subdot \T(\subdot \neg\subdot\iota (t))\vee\neg \mrm{Sent}_{\lt}(\val{t}).  
 \eeq
 
The left-to-right direction of \eqref{eq:demk2b} follows directly from $\T15$. Similarly for the right-to-left direction, if $\T\subdot \T(\subdot \neg \subdot\iota (t))$ holds. If $\neg \mrm{Sent}_{\lt}(\val{t})$, $\iota(\val{t})=\corn{0=1}$, so $\T\subdot \T\mrm{num}(\subdot\T\subdot \neg\subdot \iota(t))$. 
\end{proof}

Next, we turn to the question how $\cdplust$ and $\ckf$ are related. Corollaries \ref{cor:cdpit} and \ref{cor:cdpdt} show that some of the axioms of $\ckf$ are not provable in $\cdplus$.  However, as we shall see shortly, \emph{these are all provable} in $\cdplust$.  In fact, we shall see that $\cdplust$ and $\ckf$ are identical theories. 

\begin{lemma}\label{le_cdplust-in-ckfcs}
$\cdplust$ is a subtheory of $\ckf$.
\end{lemma}

\begin{proof}
	We verify a few key axioms, reasoning informally within $\ckf$.

\medskip
\noindent
For T2$^+$, we use T14, \eqref{eq:obs_equivalence}, and distribution of T over $\lor$.

\medskip
\noindent
For T3, since $\T\T t\ra\T\val t$ is an axiom of $\ckf$, we have $\T\T t\ra(\T\T t\lra\T\val t)$. So assume $\T\F t$. Then $\F\val t$, hence $\neg\T\val t$ and therefore $\T\val t\ra\T\T t$, hence $\T\val t\lra\T\T t$.

\medskip
\noindent
For D3, we use \textup{\ref{eq:inta}} and \eqref{eq:obs_Fdisj} -- derivable in $\ckf$ -- along with axioms T14 and T15.

\medskip
\noindent
For D5, assume $\T\T(\vphi\land\psi) \lor \T\F(\vphi\land\psi)$. By T10 and T11, this is equivalent to
	\begin{IEEEeqnarray*}{RCLCL}
		(\T\T\vphi\land\T\T\psi) &\lor& (\T\F\vphi\lor\T\F\psi).
	\end{IEEEeqnarray*}
Since we have $\T\F t\ra\F\val t$, by a series of propositional inferences we obtain
	\begin{IEEEeqnarray*}{RCLCL}
		\underbrace{((\T\T\vphi\lor\T\F\vphi)\land(\T\T\psi\lor\T\F\psi))}_{\delta_1} \ \lor \ \underbrace{((\T\T\vphi\lor\T\F\vphi)\land \F\vphi)}_{\delta_2} \ \lor \ \underbrace{((\T\T\psi\lor\T\F\psi)\land \F\psi)}_{\delta_3}.
	\end{IEEEeqnarray*}

Conversely, each combination obtainable from the conjuncts of $\delta_1$ entails $\T\T(\vphi\land\psi) \lor \T\F(\vphi\land\psi)$. As $\delta_2$ and $\delta_3$, use T16 and T4 to obtain $\T\F\vphi\land\F\vphi$, respectively $\T\F\psi\land\F\psi$, which entail the desired conclusion via T11.
\end{proof}

\begin{lemma}\label{le_ckfincd}
	$\ckf$ is a subtheory of $\cdplust$.
\end{lemma}

\begin{proof}
	The key observation is that, since $\cdplust\vdash\T\T\vphi\lor\T\F\vphi \ra (\T\vphi\lra\T\T\vphi)$, we can perform the necessary quotation and disquotation steps to prove compositionality within layers of \T. For example, to see that $\cdplust\vdash \T9$, assume $\T\T\vphi$ (or $\T\T\neg\neg\vphi$). Then $\T\T\vphi\lor\T\F\vphi$ and $\T\T\neg\vphi\lor\T\F\neg\vphi$, therefore we have both $\T\T\vphi\lra\T\vphi$ and $\T\T\neg\neg\vphi\lra\T\neg\neg\vphi$. But then $\T\T\vphi\lra\T\vphi\lra\T\neg\neg\vphi\lra\T\T\neg\neg\vphi$. The derivability of other axioms follows a similar pattern. We show some examples, reasoning informally within $\cdplust$.

\medskip
\noindent
T10: Assume $\T\T(\vphi\land\psi)$. Then $\T\T(\vphi\land\psi)\ra\T(\vphi\land\psi)$ by T3,\footnote{
	Under the definition $\D x:\lra\T\T x\lor\T\F x$. Same remark applies below. 
}
hence $\T\vphi\land\T\psi$, and therefore $\neg\F\vphi\land\neg\F\psi$. This together with $\T\T(\vphi\land\psi)$ yields, via D5, that $\vphi$ and $\psi$ are determinate, i.e., $\T\T\vphi\lor\T\F\vphi$ and $\T\T\psi\lor\T\F\psi$. Hence we conclude $\T\T\vphi\land\T\T\psi$ from $\T\vphi\land\T\psi$.

For the converse direction, $\T\T\vphi\land\T\T\psi$ entails $\T\T(\vphi\land\psi)\lor\T\F(\vphi\land\psi)$ by D5. Moreover, since $\T\T\vphi\land\T\T\psi$ also entails $\T\vphi\land\T\psi$, we get $\T(\vphi\land\psi)$, hence $\T\T(\vphi\land\psi)$ by T3.

\medskip
\noindent
T16: If $\T\T\vphi$, then $\T\T\vphi\lor\T\F\vphi$, hence $\T\vphi$ by T3.
\end{proof}

As a corollary we obtain the following

\begin{prop}\label{prop:ckfcdt}
	$\cdplust$ and $\ckf$ are identical theories. 
\end{prop}

\section{Alternative axiomatizations}\label{sec:AltAx}
	As mentioned, the theory $\ckf$ contains, along with Kripkean truth conditions in the inner layers of truth, an axiom (TDel) restricting the class of fixed-points to those that are consistent. In this section, we discuss alternative consistency axioms.

The axiom TDel is reminiscent of the schema often called T-Out: $\T\corn A\ra A$. It is also known (see \cite{can96}) that T-Out is equivalent to the consistency axiom {\sf Cons}: $\forall\vphi:\lt (\F \vphi\ra\neg\T \vphi)$. It may be asked whether TDel is equivalent to \textup{\textup{\ref{eq:inta}}} from Observation~\ref{obs_ckf-theorems}: $\T\F \vphi\ra\neg\T\T\vphi$. The fact that the latter follows from the former is straightforward: $\T\F \vphi\ra\F\vphi\ra\neg\T\vphi\ra\neg\T\T\vphi$. However, the converse employs an extra axiom stating that only sentences are truly true.
\begin{lemma}\label{le_TCons-from-(8)}
	$\ckf-\T\mrm{Del}+\textup{\textup{\ref{eq:inta}}}+\T\T t\ra\mrm{Sent}_{\lt}(t^\circ)\vdash\T\mrm{Del}$.
\end{lemma}

\begin{proof}
	We distinguish cases. If $t^\circ\notin\mrm{Sent}_{\lt}$, then $\neg\T\T t$, hence trivially $\T\T t\ra\T t$. Else, we reason by induction on the formal complexity of $t^\circ=\vphi$. If $\vphi$ is an equality, the claim follows from T7 and T8. If $\vphi\equiv\T t$, the claim follows from T14. If $\vphi\equiv\neg\T t$, we use T15 and the additional axiom: $\T\T\neg\T t\ra\T\F t\ra\neg \T\T t\ra\T\neg \T t$.
\end{proof}

Additionally, it can be observed that \textup{\ref{eq:inta}} is equivalent to $\T\F\T t\ra\F\T t$.

\begin{lemma}
	The theories $\ckf-\T\mrm{Del}+\textup{\ref{eq:inta}}$ and $\ckf-\T\mrm{Del}+\T\F t \ra \neg \T\T t$ prove the same theorems.
\end{lemma}

\begin{proof}
	To derive \textup{\ref{eq:inta}} from $\T\F\T t\ra\F\T t$, we use T15: $\T\F t\ra\T\F\T t\ra\F\T t\ra\neg\T\T t$. Conversely, the claim follows from the previous lemma observing that, for $t\in\mrm{CTerm}$, we have $\T\F t\in\mrm{Sent}_{\lt}$.
\end{proof}

Collecting these observation together, we obtain the following proposition (cf.~\cite[Lemma 15.9]{hal14}).

\begin{prop}\label{pr_Cons-equiv}
Over $\ckf\,-\,\T\mrm{Del}\, +\, \T\T t\ra\mrm{Sent}_{\lt}(t^\circ)$, the following statements are equivalent:
	\begin{enumeratei}
		\item $\T\T t\ra\T t$,
		\item $\T\F t \ra \neg \T\T t$,
		\item $\T\F\T t\ra\F\T t$.
	\end{enumeratei}
\end{prop}

Some formulations of $\kf$ (e.g. \cite{rei85}, \cite{can89}) do include the axiom $\T t\ra\mrm{Sent}(t^\circ)$. However, we decided not to include $\T\T t\ra\mrm{Sent}_{\lt}(t^\circ)$ in our official formulation of $\ckf$ in order to simplify its comparison with $\cdplust$. If $\ckf$ were defined without TDel but with, for example, \textup{\ref{eq:inta}} along with $\T\T t\ra\mrm{Sent}_{\lt}(t^\circ)$, then the equivalence stated in Proposition~\ref{prop:ckfcdt} would need to be reformulated as follows:

\begin{lemma}
	The theories $\ckf':=\ckf - \T\mrm{Del} + \textup{\ref{eq:inta}} + \T\T t\ra\mrm{Sent}_{\lt}(t^\circ)$ and $\cdplust + \D t \ra \mrm{Sent}_{\lt}(t^\circ)$ are identical.\footnote{
		The same would hold for $\cdplust + \T\T t \ra \mrm{Sent}_{\lt}(t^\circ)$.
	}
\end{lemma}

\begin{proof}
	For the inclusion of $\cdplust + \D t \ra \mrm{Sent}_{\lt}(t)$ into $\ckf'$, the crucial axioms are T3 and $\D t \ra \mrm{Sent}_{\lt}(t^\circ)$. As for the latter, $\T\T t\ra\mrm{Sent}_{\lt}(t^\circ)$ is just an axiom of $\ckf'$. As for $\T\F t\ra\mrm{Sent}_{\lt}(t^\circ)$, it follows from $\T\T\ng t\ra\mrm{Sent}_{\lt}(\ng t^\circ)$, hence $\mrm{Sent}_{\lt}(t^\circ)$.
	
	For T3, we distinguish two cases. If $t^\circ\notin\mrm{Sent}_{\lt}$, then $\ng t^\circ\notin\mrm{Sent}_{\lt}$, hence $\T\T t\ra\bot\land\T\F t\ra\bot$, hence trivially $\T\T t\lor\T\F t\ra(\T\T t\lra\T t^\circ)$. If $t\in\mrm{Sent}(t)$, by Lemma~\ref{le_TCons-from-(8)} we have $\T\T t\ra\T t^\circ$, hence $\T\T t\lra(\T t^\circ\lra\T\T t)$. For the second disjunct, using Proposition~\ref{pr_Cons-equiv} we have $\T\F t\ra\F t^\circ\ra\neg\T t^\circ\ra (\T t^\circ\ra\T\T t)$.
	
	For the converse inclusion of $\ckf'$ into $\cdplust + \D t \ra \mrm{Sent}_{\lt}(t)$, the crucial cases are \textup{\ref{eq:inta}} and $\T\T t\ra\mrm{Sent}_{\lt}(t^\circ)$. The latter follows from $\D t \ra \mrm{Sent}_{\lt}(t)$. The former can be derived thus: $\T\F t\ra\F t\ra\neg\T t\ra\neg\T\T t$.
\end{proof}

\section{Complete, symmetric, and mixed fixed-points}
In this section, we verify whether the results from previous sections carry over if one focuses on different classes of fixed-points of $\mathscr{K}$. We provide a positive answer for the class of complete fixed-points as well as for the class of consistent \emph{or} complete fixed-points. Specifically, given a complete fixed-point $X$ of $\mathscr{K}$, the structure $(\nat,\mbb T_X)$ can be shown to be a model of the determinateness axioms of $\cdplus$ with $\D$ defined in terms of $\T$. Moreover, the resulting theory can be shown to be identical to a variant of $\ckf$. Similarly for the class of fixed-points which are either consistent or complete. The question whether similar results are available for mixed fixed-points, i.e., fixed-points which are neither consistent nor complete, will be left open.

\subsection{Complete models}

Since the structure of the arguments is very similar, we limit ourselves to highlighting the necessary modifications.

\begin{dfn}[$\cdplustcp$]
    Let $\cdplustcp$ be the $\lt$-theory whose axioms are those of $\cdplus$, but where $\D$ defined in terms of $\T$ as
	\begin{IEEEeqnarray}{L}
		\D x:\lra \neg\T\T x\lor\neg\T\F x. \label{eq_Ddf-comp}
	\end{IEEEeqnarray}
\end{dfn}

It can be shown that, given a complete fixed-point $X=\mathscr{K}(X)$, the structure $(\mbb N, \mbb T_X)$, where $\mbb T_X:=\{A\in \lt\mid(\mbb N, X)\vDash A\}$, is a model of $\cdplustcp$. The argument follows the blueprint of Proposition \ref{prop_cons-mod}. We first observe that, since a complete fixed-point $X$ of $\mathscr K$ is such that $(\nat, X)\vDash\kfcom$, and since $\kfcom$ derives the schema $A\ra\T\corn A$, we have

\begin{fact}\label{fact_aux-comp}
	For any complete fixed-point $X$, for any $A$, $(\nat, X)\vDash A\ra \T\corn A$.
\end{fact}

Just as Fact \ref{fact_auxiliary} was used in the proof of Proposition~\ref{prop_cons-mod}, Fact~\ref{fact_aux-comp} will play a similar role in the proof of the following

\begin{prop}\label{prop_comp-mod}
	For a complete fixed-point $X$, $(\nat,\mbb T_X)\vDash \cdplustcp$.
\end{prop}
\begin{proof}[Proof Sketch]
    By induction on the length of proofs in \cdplustcp. We verify T3 and D3, whose arguments are symmetric to those in the proof of Proposition \ref{prop_cons-mod}.
    
\medskip
\noindent
T3: Assume $(\nat,\mbb{T}_{X})\vDash\neg\T\T\vphi \lor \neg\T\F\vphi$, which is the case iff $\{\vphi,\neg\vphi\}\not\subseteq X$. To show $(\nat, \mathbb{T}_X)\vDash\T\T\vphi\ra\T\vphi$, let $\vphi=\corn A$ and assume $(\nat, \mathbb{T}_{X})\vDash\T\T\vphi$, which is equivalent to $\vphi\in X$. Towards a contradiction, suppose $(\nat, \mbb T_X)\not\vDash\T\vphi$, hence $(\nat, X)\not\vDash A$. By Fact \ref{fact_aux-comp}, we get $(\nat, X)\vDash\F\vphi$ iff $\neg\vphi\in X$, contradicting our assumption.

Conversely, to show $(\nat, \mathbb{T}_X)\vDash\T\vphi\ra\T\T\vphi$, we reason as follows without using the assumption on the determinateness of $\vphi$:
	\begin{IEEEeqnarray*}{LL}
		(\nat, \mathbb{T}_{X})\vDash\T\vphi, \hspace{3mm} & \text{ iff }
		\\
		(\nat,{X})\vDash A, & \text{ \textit{hence}, by Fact \ref{fact_aux-comp}}
        \\
        (\nat, X)\vDash \T\vphi, & \text{ iff}
        \\
        (\nat, \mathbb{T}_{X})\vDash\T\T\vphi.
	\end{IEEEeqnarray*} 

\medskip
\noindent
D3: We need to show that $(\mbb{N}, \mbb T_X)\vDash\neg\T\T\vphi\lor\neg\T\F\vphi$ is equivalent to
	\begin{IEEEeqnarray*}{RCL}
		(\mbb{N}, \mbb T_X)\vDash\neg\T\T \big(\neg\T\T\vphi\lor\neg\T\F\vphi\big) \lor 
		\neg\T\F \big(\neg\T\T\vphi\lor\neg\T\F\vphi\big).
	\end{IEEEeqnarray*}
The latter is equivalent to
	\begin{IEEEeqnarray*}{RCL}
		\underbrace{(\mbb{N}, X)\vDash\neg\T\big(\neg\T\T\vphi\lor\neg\T\F\vphi\big)}_{d_1}
		\ &\text{ or }& \ 
		\underbrace{(\mbb{N}, X)\vDash\neg\F\big(\neg\T\T\vphi\lor\neg\T\F\vphi\big)}_{d_2}.
	\end{IEEEeqnarray*}

Since $d_1$ and $d_2$ are equivalent to, respectively, $(\mbb{N}, X)\vDash\neg\T\vphi\land\neg\F\vphi$ and $(\mbb{N}, X)\vDash\neg\F\vphi\lor\neg\T\vphi$, by completeness of $X$ it can be observed that their disjunction is equivalent to $d_2$, which is equivalent to $(\mbb{N}, \mbb T_X)\vDash\neg\T\T\vphi\lor\neg\T\F\vphi$, as required.
\end{proof}

A collection of principles inspired by the models $(\mbb N, \mbb T_X)$ for $X$ a complete fixed-point can be obtained by modifying $\ckf$ in the obvious way, that is, by replacing the axiom expressing consistency with one expressing completeness, leaving the remaining axioms characterising the closure of {\sf KF} unmodified:
\begin{dfn}[$\ckfcp$]
    $\ckfcp$ is the system obtained from $\ckf$ by replacing TDel with 
    \begin{IEEEeqnarray}{L}
		\T t\ra \T\T t. \ztag{TRep}
	\end{IEEEeqnarray}
\end{dfn}
Axiom TRep readily yields completeness $\T\T t\lor\T\F t$ via $\neg\T\T t\ra\neg\T\val t\ra\F\val t\ra\T\F t$.

\begin{prop}
    \cdplustcp and \ckfcp are identical theories.
\end{prop}

\begin{proof}[Proof Sketch]
    As in the proof of Lemma~\ref{le_ckfincd}, the key observation to show that \ckfcp is a subtheory of \cdplustcp is that $\cdplustcp\vdash\neg\T\T\vphi\lor\neg\T\F\vphi\ra(\T\vphi\lra\T\T\vphi)$. We can derive the counterpositive of each axiom of $\ckfcp$ by performing the necessary quotation and disquotation step. For example, for T10 we assume $\neg\T\T(\vphi\land\psi)$, which implies $\neg\T(\vphi\land\psi)$, hence $\neg\T\vphi\lor\neg\T\psi$. Via D5, at least one between $\vphi$ and $\psi$ is determinate, hence we conclude $\neg\T\T\vphi\lor\neg\T\T\psi$. In a similar way we can derive the counterpositive of TRep: if $\neg\T\T t$, then $\T t\lra \T\T t$, hence $\neg\T t$.

    Conversely, to show that \cdplustcp is a subtheory of \ckfcp, the reasoning is similar to the proof for Lemma \ref{le_cdplust-in-ckfcs}. We only verify T3: Since $\T t\ra\T\T t$, we have immediately $\neg\T\T t\ra(\T t\lra\T\T t)$; for the other disjunct, $\neg\T\F t\ra\neg\F t\ra\T t\ra(\T\T\ra\T t)$.
\end{proof}

In fact, the well-known duality between consistency and completeness 
is preserved in the present setting, in that $\ckf$ and $\ckfcp$ are mutually $\lnat$-interpretable via Cantini's dual translation, mapping $\T$ to $\neg\F$ \cite{can89}.
More precisely, let $c$ be a map of $\lt$ into itself preserving the arithmetical vocabulary, commuting with logical operations in the usual way, and mapping $\T x$ to $\neg\F x$.

\begin{prop}
	$\ckf$ and $\ckfcp$ are mutually $\lnat$-interpretable via $c$.
\end{prop}

\begin{proof}[Proof Sketch]
	For R1 and T1-T15, it suffices to observe that, within both $\ckf$ and $\ckfcp$, their instances $A$ are self-dual, in the sense that $A\lra A^c$. For example, for T15 we have (tacitly using the fact that $t\in\mrm{Cterm}$)
	\begin{IEEEeqnarray*}{RCL+L}
		(\T\F\T t \lra \T\F t)^c &\text{ iff }& \neg\F\neg\F\neg\neg\T\neg t \lra \neg\T\neg\neg\F\neg t & \text{def of $c$}
		\\
		&\text{ iff }& \neg\F\neg\F\T\neg t \lra \neg\T\F\neg t & \text{T9}
		\\
		&\text{ iff }& \neg\T\F\T\neg t \lra \neg\T\F\neg t & \text{def of $\F$, T9}
	\end{IEEEeqnarray*}
The last line is a counterpositive instance of T15.

Similarly for TDel and TRep, just note that $\text{TDel}^c=\neg\F\neg\F t\ra\neg\F t\lra(\neg\T\F t\ra\neg\F t)$, which is the counterpositive of TRep, and $\text{TRep}^c=\neg\F t\ra\neg\F\neg\F t\lra(\neg\F t\ra\neg\T\F t)$, which is the counterpositive of TDel.
\end{proof}

The duality between $\msf{KF}+\msf{CONS}$ and $\msf{KF}+\msf{COMP}$ can be lifted to intertranslatability a.k.a.~synonymy (see \cite{nic21dua}); the same holds for $\ckf$ and $\ckfcp$.

\begin{corollary}
$\ckf$ and $\ckfcp$ are synonymous. 
\end{corollary}
\noindent In particular, this means that the logical strength of $\ckfcp$, too, coincide with that of $\cdplus$.

\subsection{Symmetric models}

Combining the results on consistent and complete fixed-points, it can also be shown that the definition of $\D$ can be adapted to obtain that the class of \textit{symmetric} (i.e., consistent or complete) fixed-points of $\mathscr{K}$ satisfies the axioms of $\cdplust$ with $\D$ defined disjunctively as follows:

\begin{dfn}[$\cdplustsym$]
    Let $\cdplustsym$ be the $\lt$-theory whose axioms are those of $\cdplus$, but where $\D$ defined in terms of $\T$ as
	\begin{IEEEeqnarray}{L}
		\D x:\lra (\T\T x\lor\T\F x) \land (\neg\T\T x\lor\neg\T\F x). \label{eq_Ddf-sym}
	\end{IEEEeqnarray}
\end{dfn}

%

Using Facts \ref{fact_auxiliary} and \ref{fact_aux-comp}, it can then be shown that symmetric fixed-point $X$ can be used to obtain models of \cdplustsym:

\begin{prop}
	For any symmetric fixed-point $X$, $(\mbb N, \mbb T_X)\vDash\cdplustsym$.
\end{prop}

\begin{proof}
	By induction on the length of proofs in $\cdplustsym$. For T3, one assumes $(\nat,\mbb{T}_{X})\vDash\T\T\vphi \lor \T\F\vphi\land(\neg\T\T\vphi \lor \neg\T\F\vphi)$ and may then reason by cases, depending on whether $X$ is consistent or complete, following the arguments for Proposition~\ref{prop_cons-mod} and Proposition\ref{prop_comp-mod}.
\end{proof} 

Accordingly, a corresponding {\sf CKF} system is obtained by a disjunction of TDel and TRep:

\begin{dfn}[$\ckfsym$]
    $\ckfsym$ is the system obtained from $\ckf$ by replacing TDel with 
    \begin{IEEEeqnarray}{L}
		(\T\T t\ra \T\val t) \lor (\T\val s\ra \T\T s). \ztag{TSym}
	\end{IEEEeqnarray}
\end{dfn}

To see that axiom TSym yields \textit{consistency-or-completeness} $(\T\T t\land\T\F t)\ra(\T\T s\lor\T\F s)$, assume $\T\T t\land \T\F t$. If $(\T\T t\ra \T\val t)$, then $\T t\land \F t$, which is impossible, hence $(\T\val s\ra \T\T s)$ and therefore $\neg\T\T s\lor\neg\T\F s$.

It can also be observed that the translations $\tau$ and $\iota$ defined in the proof of Proposition \ref{pr_ckf-ctcon-cdplus} yield mutual interpretability of $\ckfsym$ with $\msf{CT}[\msf{KF}+(\msf{CONS\lor COMP})]$, hence we obtain the following equivalence:

\begin{corollary}
	The theories $\cdplus$, $\ckf$, $\ckfcp$, and $\ckfsym$ have the same arithmetical consequences as the system $\msf{RT}_{<\varepsilon_{\varepsilon_0}}$ of ramified truth up to $\varepsilon_{\varepsilon_0}$.
\end{corollary}

\subsection{Mixed Models.} Mixed fixed-point models are those in which the truth predicate can feature \emph{both} gaps and gluts. Can \textit{mixed} (i.e., neither consistent nor complete) fixed-points model $\cdplust$ under a suitable definition of $\D$? We leave this question open. However, we observe that mixed fixed-points are not models of $\cdplust$ under any of the definitions of $\D$ considered above.
\begin{prop}\label{prop:mixed}
	Let $X$ be a mixed fixed-point and let $\cdplust\star$ range over the theories $\cdplust, \cdplustcp, \cdplustsym$. Then $(\nat,\mbb{T}_{X})\not\vDash\cdplust\star$.
\end{prop}
\begin{proof}
	In light of Proposition~\ref{prop_cons-mod}, we show that $(\nat,\mbb T_X)\not\vDash\T3\lor\D5$,\footnote{%
		This is of course redundant, but it clarifies the reasons why each of the axiom is not satisfied.
	}
which are the two axioms where Facts~\ref{fact_auxiliary} and \ref{fact_aux-comp} played a crucial role. 

To show $(\nat,\mbb{T}_{X})\not\vDash\T3$, let $\lambda$ and $\tau$ be such that $\lambda\land\neg\lambda\in X$ and $\tau\lor\neg\tau\notin X$, and assume moreover that $\mrm{PAT}\vdash\lambda\lra\neg\T\lambda$ and $\mrm{PAT}\vdash\tau\lra\T\tau$. Since $\neg(\tau\lor\lambda)\notin X$, it can be observed that the following jointly hold
	\begin{IEEEeqnarray*}{LL}
		(\nat,{X})\vDash\neg\F(\tau\lor\lambda),
		\\
		(\nat,{X})\vDash\T(\tau\lor\lambda),
		\\
		(\nat,{X})\vDash\neg(\tau\lor\lambda).
	\end{IEEEeqnarray*}
We derive, for $\D$ defined by \eqref{eq_Ddf-sym},
	\begin{IEEEeqnarray*}{L}
		(\nat,\mbb{T}_{X})\vDash\D(\tau\lor\lambda)\land\T\T(\tau\lor\lambda)\land\neg\T(\tau\lor\lambda),
	\end{IEEEeqnarray*}
hence $(\nat,\mbb{T}_{X})\not\vDash\T3$.

To show $(\nat,\mbb{T}_{X})\not\vDash\D5$, let similarly $\vphi$ be such that $\vphi\land\neg\vphi\in X$, and let $\psi$ be such that $\psi\lor\neg\psi\notin X$. This entails
	\begin{IEEEeqnarray*}{LL}
		(\nat,{X})\vDash\F(\vphi\land\psi) \land \neg\T(\vphi\land\psi), & \text{ iff}
		\\
		(\nat,\mbb{T}_{X})\vDash\T\F(\vphi\land\psi)\land\neg\T\T(\vphi\land\psi), \hspace{1mm} & \text{ hence}
		\\
		(\nat,\mbb{T}_{X})\vDash\D(\vphi\land\psi).
	\end{IEEEeqnarray*}
with $\D$ defined as per \eqref{eq_Ddf-sym}. However, $(\nat,\mbb{T}_{X})\not\vDash\D\vphi\lor\D\tau$ for any of the above definitions of $\D$, since $(\nat,\mbb{T}_{X})\vDash(\T\T\vphi\land\T\F\vphi)\land(\neg\T\T\psi\land\neg\T\F\psi)$, hence $(\nat,\mbb{T}_{X})\not\vDash\D5$.
\end{proof}

The reason why the class of symmetric, but not that of mixed, fixed-points is suitable for modeling $\cdplust$ can be explained as follows. The former, but not the latter, features a specific interplay between the notions of \textit{determinate} and \textit{having a classical semantic value}. Within both the class of consistent and the class of complete fixed-points, we can single out an intended interpretation, where being determinate can be defined as having a classical semantic value. These intended models are the least fixed-point, and the largest fixed-point, where the set of sentences with a classical value are those which are grounded in Kripke's sense. Proposition \ref{prop:mixed} clarifies  why this is not possible in mixed models: in each of them, Boolean combinations of gluts and a gaps result in sentences which are \textit{strictly} true or \textit{strictly} false.\footnote{This is essentially the same reason why the version of $\kf$ in \cite{fef91} cannot satisfy the $\D$-axioms of $\cd$ and $\cdplus$ when
\[
\D x :\lra (\T x\vee \F x) \land \neg (\T x\land \F x).
\]}

\section{Assessment}\label{sec:ass}

The theories introduced in this work have unique features that place them among the most promising theories of truth available in the literature. In this section, we elaborate on some of these features, comparing variants of $\ckf$ with related approaches to truth and determinateness along the key dimensions outlined in the introductory section.

The generalizing function of truth, we argued, requires a strongly classical and fully compositional theory of truth. In our theories, such a function is realized in virtue of axioms $\T4$-$\T6$.
For example, $\T4$ is sufficient to exclude the existence of sentences that are both true and false and of sentences that are neither true nor false. In this sense, all variants of both $\cd$ and $\ckf$ are \textit{classical theories of classical truth}. To state this formally, recall that the \emph{internal theory} of a theory of truth $S$ be defined as 
	\[
		\mathscr{I}(S)=\{ A \sth S \vdash \T\corn{A}\}.
	\]
For every variant $\msf{CKF}$ and $\cd$, the \emph{logic of} their internal theories is \emph{classical}: (i) all (universal closures of) classical logical axioms in $\lt$ are true, and (ii) all classical logical inferences preserve truth, hence, (iii) all (universal closures of) theorems of classical logic in $\lt$ are true.\footnote{For this induction to hold, it is important that the induction schema of $\ckf$ and variants is extended to $\T$.}

By contrast, systems such as $\kf$ or $\msf{DT}$ are \textit{classical theories of nonclassical truth}: they are formulated in classical logic, yet the logics governing their internal theories is nonclassical. As noted in the introduction, this compromises the generalizing power of its truth predicate. 
Advocates of such theories -- e.g.~\cite{rei86} -- stress that their internal theories, despite not obeying the laws of classical logic, enjoy other truth theoretic virtues. For instance, the internal theory of $\kf+\msf{CONS}$ is not only closed under Strong Kleene Logic, but it is \textit{fully disquotational} too: $A\in\mathscr{I}(\kf+\msf{CONS})$ iff $\T\corn A\in\mathscr{I}(\kf+\msf{CONS})$. Full disquotation is often taken to be crucial by truth theorists (see, e.g., \cite{fie08}). This property provably fails for the inner theory of any classical theory of classical truth that admits a standard model.

However, our results show that, in addition to the internal theory, another notion plays prominent theoretical role in this context. Define the \emph{deep theory} of $S$ as
	\[
 	   \mathscr{D}(S)=\{ A \sth S\vdash \T\corn{\T\corn{A}}\}.
	\]
In theories of classical determinate truth such as the ones studied in this paper, the internal and deep theories \emph{provably differ}.\footnote{
	For instance, in $\ckf$, $\lambda\vee\neg\lambda$ can be used to separate the two. Dually, one can use $\lambda\land\neg\lambda$ to separate the inner and the internal theory of $\ckfcp$, since the latter theory derives $\T\T(\lambda\land\neg\lambda)$. 
}
One of the main virtues of the theories we propose is that the logics of their deep theories can be associated with well-known logics admitting a transparent truth predicate. For $\ckf$, the logic of its deep theory amounts to the familiar Strong Kleene logic with a transparent truth predicate. Thus, although not every classical axiom will be in $\mathscr{I}(\ckf)$, whatever \textit{is} inside it can be closed under the relevant non‑classical rules of inferences and under iterations of $\T$. Analogously, the logic of $\mathscr{D}(\ckfcp)$ corresponds to the Logic of Paradox \cite{priLogic}, and the logic of $\mathscr{I}(\ckfsym)$ to Symmetric Strong-Kleene \cite{scoCombinators, bla02}. 

As anticipated, our theories also provide clear semantic rules for the analysis of the language with type-free truth. To see this, we note that we can uniformly define a `semantic' truth predicate $\T_\mrm{sem}$ as $\T_\mrm{sem} x:\lra \T\T x$. For such a predicate, the theories can prove universally quantified laws corresponding to the Strong Kleene truth conditions for $\lt$.\footnote{Essentially, the compositional axioms of $\kf$.} In addition, the theories prove unrestricted positive and negative truth ascriptions for $\T_\mrm{sem}$. In formalizing semantic rules for $\lt$, our theories also provide definite information on the space of `models', or extensions of $\T_\mrm{sem}$, that are admissible. While $\ckf$ only allows consistent interpretations of $\T_\mrm{sem}$, $\ckfcp$ forces inconsistent but complete interpretations. 

This is in stark contrast with $\cd$ and its variants,
as they do not prove the Strong Kleene conditions for $\T_\mrm{sem}$ 
(cf.~Corollaries \ref{cor:cdpdt} and \ref{cor:cdpit}). This means that, in such theories, truth behaves transparently and according to logical principles only on a restricted fragment of the language, namely on determinate sentences. 
In addition, $\cd$ and its variants do not impose clear conditions of admissible interpretations of the true and determinate sentences -- the analogue of our $\T_\mrm{sem}$. As the results in Fujimoto and Halbach show \cite[Thm.~7.11 and Thm.~8.1]{fuha23}, the theories are compatible with consistent or inconsistent interpretation of the true and determinate sentences. In our approach, there is a clear choice to be made depending on one's chosen definition of determinateness: while the axioms of $\cdplus$ formulated by means of the definition of $\D$ as $\T\T x\vee \T\F x$ result in the theory $\ckf$, the dual definition as $\neg \T\T x\land \neg \T\F x$ yields the theory $\ckfcp$ whose deep theory is inconsistent.

Combining these observations together, we see that variants of $\ckf$ (i) are fully compatible with the generalizing function of truth -- unlike classical theories of non‑classical truth, and (ii) capture a transparent and well-behaved notion of truth inside their deep theory -- unlike $\cd$ and its variants.


Another fundamental feature of the theories introduced in this work is that the notion of determinateness, just like what happens in well-known formal approaches to truth, is  \emph{defined} in terms of truth. Feferman \cite{fef08}, for example, despite assigning priority to the axioms for determinateness over those for truth, defines a sentence to be determinate iff it is true or false (and not both). A similar case for a determinateness predicate defined in this way can be made for $\kf$ -- see especially \cite{rei86} and \cite{fef91}. 

Feferman's definition of determinateness is well suited for theories like $\mathsf{DT}$ and $\kf$, which employ a self-applicable but nonclassical truth predicate. As Fujimoto and Halbach rightly point out, however, the definition $\D x:\lra \T x\vee \F x$ is not appropriate for theories of a thoroughly classical conception of truth, such as $\cd$ and $\cdplus$, or classical closures of Kripke-Feferman truth developed in this paper. This is why $\cd$ and $\cdplus$ treat determinateness as a primitive.

By contrast, in our theories a strongly classical and compositional truth predicate co-exists with a {defined} determinateness predicate. In particular, our results show that the desiderata imposed to the notion of truth and determinateness by theories such as $\cdplus$ can in fact be realized by theories based on a defined determinateness predicate. While the definition $\D x:\lra \T x\vee \F x$ is unsuitable in this context, the alternative $\D x:\lra \T \T x \vee \T \F x$ is just right to license the principles of $\cdplus$, and more generally to meet the core desiderata for a thoroughly classical, self-applicable conception of truth.

One might object that our definition of $\D$ seems more artificial than Feferman's, which rests on the natural thought that being determinate just means `having a determinate (classical) truth-value' in a paracomplete or paraconsistent model. But the definition we propose is in fact a rather natural incarnation of this standard notion. In particular, it puts $\ckf$'s determinateness in continuity with the notion of determinateness available within $\kfcon$. 
Since in $\kfcon$ determinateness is defined as $\T x \lor \F x$, extending the theory to its classical closure naturally requires introducing an additional layer of truth into the definition of $\D$. If one endorses Feferman's extension of determinateness as given in various manifestations of the Kripke-Feferman theory, its extensions will remain unchanged in our theories. What changes is the generalizing power afforded by the classical truth-theoretic layer. More precisely, the $\nat$-categoricity of our theories (Proposition \ref{pr_Ncat}) tells us that any $\omega$-model of our theory features a standard, Kripke-Feferman notion of determinateness. And in each such model compositional axioms and classical logic are fully satisfied in a strong sense because it amounts to a classical closure of a Kripkean fixed point. In particular, if we consider the classical closure of the minimal fixed point of Kripke's theory of truth, the extension of its determinateness predicate is just the set of grounded sentences of $\lt$.\footnote{This makes fully transparent a link adumbrated by Fujimoto and Halbach, when they ``call the sentences in $\mathfrak{D}$ \textit{determinate}. If it were not for the additional predicate $\D$, they would be the sentences that are grounded in Kripke's \cite{kri75} sense (with some qualifications)'' \cite[p.~222]{fuha23}. In our setting, this connection is made completely explicit without the need to reinterpret sentences featuring the primitive determinateness predicate.}

We can reformulate the above observations without reference to Kripkean semantics. The Kripke-Feferman notion of determinateness amounts to being true or false. Its extension remains unchanged in (the consistent variant) of our theories, but its definition takes into account the strong classicality of our truth predicate required by its generalizing role. Determinateness is now defined as being \emph{semantically true or false}, where `semantically true' can be explicitly defined as the predicate $\T_\mrm{sem}$ just introduced. 


Taken together, the above observations point to a novel approach to truth and determinateness, yielding a family of axiomatic theories with distinctive features. Known theories that define determinateness in terms of truth and falsity -- such as $\kf$ or $\msf{DT}$ -- employ a compositional and self-applicable truth predicate, but one that, due to its nonclassical nature, does not perform well on generalizations such as the ones required in blind deductions.\footnote{See again \cite{fuj22}.} Conversely, theories of classical truth like the Tarskian theory $\msf{CT}$, variants of $\cd$, or Friedman and Sheard's $\msf{FS}$, feature a strongly
compositional truth predicate, but one that is not suitable to capture the Kripke-Feferman notion of determinateness. In addition, for $\msf{CT}$ and $\msf{FS}$, {no equally satisfactory notion of determinateness is likely to be found.}

However, the classical closures of $\kf$ show that it is possible
to combine, in a single framework, both kinds of virtues: a classical, strongly compositional notion of truth that supports blind inferences, and the possibility of defining a class of determinate sentences satisfying desirable principles.\footnote{
    As already mentioned in the introduction, a thorough philosophical assessment the conceptual import of $\ckf$ will be carried out in future work. Our focus in this paper was on its formal properties.
}

\section{Extensions and open questions}

To conclude, we list some open questions and potential lines of research stemming from our work.

\begin{enumerate}[label=(\arabic*)]
    \item Each of the introduced variants of $\ckf$ restrict, in different ways, the class of Kripkean fixed-points. For example, the axiom \text{TDel} in $\ckf$ restricts the fixed-points to the consistent ones. Let $\msf{CKF}:=\ckf-\text{TDel}$. Then:
    \begin{enumeratei}
        \item Is there a definition of $\D$ that ensures the identity of $\msf{CKF}$ and the appropriate reformulation of $\cdplust$ in terms of such a definition?
        \item Would this definition of $\D$ ensure the soundness of $\cdplust$ with respect to mixed models?
    \end{enumeratei}
\end{enumerate}

Moreover, one can ask whether there are natural principles that can be used to expand $\ckf$, and how strong the resulting theories become. We mention two possible such expansions.

\begin{enumerate}
\setcounter{enumi}{1}
    \item Given the semantics of $\ckf$, one might treat the least fixed-point of $\mathscr{K}$ as the intended interpretation for its deep theory. One can then axiomatize this conception by adding a minimality schema to $\ckf$. So let $\msf{CKF}_\mu$ be the theory obtained by expanding $\ckf$ with the schema
    		\[
    			K(\vphi(x)) \ra \forall x(\T\T x\ra\vphi(x)),
    		\]
    	where, for $\vphi\in\mathcal L_\T$, $K(\vphi(x))$ expresses that $\vphi(x)$ is a closed point of $\mathscr{K}$. Then:
    
   \medskip	
    \begin{enumeratei}
    \setcounter{enumii}{2}
    		\item  What is the proof-theoretic strength of $\msf{CKF}_\mu$?
    \end{enumeratei}
    \medskip
    
    \item One may also consider extensions of $\ckf$ by means of refection principles. It is well-known that the theory of truth enables one to directly express soundness extensions of theories in the form of \emph{Global Reflection Principles} for a theory $S$: for any $\vphi$, if $\vphi$ is provable in $S$, then $\T\vphi$. It is also well-known that the most prominent type-free theories of truth don't sit well with Global Reflection.\footnote{
    For instance, the addition of the Global Reflection to $\kf$ results in internal inconsistency, while $\kfcon$ and $\mathsf{FS}$ are outright inconsistent with Global Reflection. Relatedly, $\ckf$ and $\cd$ are inconsistent with 
the $\omega$-iterated Global Reflection Principle.
    } 
However, for such theories Global Reflection may not be the right soundness extension to focus on. Reinhardt \cite{rei86} proposed to consider instead a partial reflection principle for \emph{provably true} sentences.
    A similar strategy is available in $\ckf$ and its variants, if one looks at theorems under two layers of truth. The principle
\beq\label{eq:grptt}
(\forall \vphi \colon \lt)(\mathrm{Prov}_{\ckf}(\T\T\vphi)\ra \T\T\vphi)
\eeq
is sound with respect to the semantics for $\ckf$ and can be safely added (and iterated) over our theories. It can then be asked:

\medskip
	\begin{enumeratei}
	\setcounter{enumii}{3}
		\item What is the strength of the reflection principle \eqref{eq:grptt} and iterations thereof?\footnote{A proof-theoretic analysis of reflection principles of this kind is currently being conducted in joint work of [omissis] and the first-named author of this paper.}
	\end{enumeratei}
    \item Finally, for technical interest, it remains to be explored whether an analogue of our results is available for the theory $\cd_{\T}$---as well as its variant $\cd_{\T}[{\sf COMP}]$. Specifically:
    \begin{enumeratei}
	\setcounter{enumii}{4}
        \item Does $\cd_{\T}$ have a natural semantics validating the principle T2, i.e., $\D t\ra\T\D t$, but not necessarily its converse?
        \item Can the semantic construction be axiomatized in such a way that the resulting theory is identical to $\cd_{\T}$?
        \item What is the proof theoretic strength of the resulting systems?
    \end{enumeratei}
\end{enumerate}

\subsection*{Acknowledgements} Luca Castaldo's work was supported by the Italian Ministry of Education, University and Research, PRIN 2022 program ``HERB - Human Explanation of Robotic Behaviour'' prot. 20224X95JC, and the National Science Centre (NCN), Poland, MAESTRO grant 2019/34/A/HS1/00399 “Epistemic and Semantic Commitments of Foundational Theories.” Carlo Nicolai's work was supported by the UK Arts and Humanities Research Council (project AH/V015516/1) and the European Union’s Horizon research and innovation programme within the projects PLEXUS (Grant agreement No 101086295). Both authors would like to thank Kentaro Fujimoto and Volker Halbach for discussing with the authors the material contained in the paper, and the audiences of the KCL-Bristol Logic Meeting, the FOMTIC conference at University of Warsaw, the PhilMath Seminar at the IHPST Paris 1, the XIII Workshop in Philosophical Logic at CONICET, Buenos Aires.




\bibliographystyle{alpha}
\bibliography{Aux/cd_bib}
 \end{document}